\documentclass[a4paper,english, 11pt]{article}
\usepackage[a4paper, margin=1.3in
]{geometry}

\usepackage[utf8]{inputenc}
\usepackage[T1]{fontenc}
\usepackage{lmodern}
\usepackage[english]{babel}
\usepackage{hyperref}

\usepackage{amsmath}
\usepackage{amsthm}
\usepackage{amsfonts}
\usepackage{amssymb}
\usepackage{mathtools}
\usepackage{siunitx}
\usepackage[linesnumbered, ruled, noend]{algorithm2e}
\usepackage{graphicx}
\usepackage{subcaption}
\usepackage{longtable}
\usepackage{multirow}
\usepackage{enumerate}
\usepackage{xcolor}
\usepackage{xspace}
\usepackage{tikz}
\usepackage[font=small,labelfont=bf]{caption}
\usepackage{array}
\usepackage{setspace}

\usetikzlibrary{fit}
\usetikzlibrary{shapes}
\usetikzlibrary{shapes.geometric}
\usetikzlibrary{backgrounds}
\usetikzlibrary{arrows}
\usetikzlibrary{arrows.meta}
\tikzset{
	-Latex,auto,node distance =0.5 cm and 0.5 cm,semithick,
	square/.style={regular polygon,regular polygon sides=4}
}

\usetikzlibrary{decorations}
\usetikzlibrary{decorations.text}
\usetikzlibrary{topaths}
\usetikzlibrary{patterns}
\usetikzlibrary{calc}
\usetikzlibrary{intersections}
\usetikzlibrary{through}
\usetikzlibrary{spy}
\usetikzlibrary{matrix}
\usetikzlibrary{lindenmayersystems}
\usetikzlibrary{external}
\usetikzlibrary{positioning}
\usetikzlibrary{3d}

\theoremstyle{plain}
\newtheorem{theorem}{Theorem}[section]
\newtheorem{lemma}[theorem]{Lemma}
\newtheorem{corollary}[theorem]{Corollary}

\newtheorem{definition}{Definition}[section]

\usepackage[pagewise]{lineno}

\usepackage[
backend=biber,
style=apa,
urldate=short,
]{biblatex}

\AtNextBibliography{\small}

\DefineBibliographyStrings{english}{%
  urlseen = {Accessed on}, 
}

\DeclareFieldFormat{urldate}{(\bibstring{urlseen}\space#1)}
\DeclareFieldFormat{url}{\url{#1}}

\renewbibmacro*{url+urldate}{%
  \usebibmacro{url}
  \usebibmacro{urldate}%
}

\usepackage{csquotes}
\usepackage{authblk}

\makeatletter
\renewcommand{\@fnsymbol}[1]{\ifcase#1\or i\or ii\or iii\or iv\or v\or vi\or vii\or viii\or ix\or x\else\@ctrerr\fi}
\makeatother

\newcommand{\python}{{\sc python}\xspace}

\newcommand{\gurobi}{{\sc gurobi}\xspace}
\newcommand{\gurobiVersion}[2]{{\sc gurobi}~\oldstylenums{#1.#2}\xspace}

\newcommand{\cplexVersion}[2]{{\sc cplex}~\oldstylenums{#1.#2}\xspace}

\newcommand{\scipjack}{{\sc \mbox{scip-Jack}}\xspace}
\newcommand{\scip}{{\sc scip}\xspace}
\newcommand{\scipVersion}[3]{{\sc scip}~\oldstylenums{#1.#2.#3}\xspace}

\newcommand{\setting}[1]{\texttt{#1}}

\newcommand{\R}{\mathbb{R}}

\newcommand{\arcVariable}{\ensuremath{x}\xspace}

\newcommand{\arcVar}{\ensuremath{x}\xspace}
\newcommand{\nodeVar}{\ensuremath{y}\xspace}
\newcommand{\bilinVar}{\ensuremath{m}\xspace}
\newcommand{\interfTotVar}{\ensuremath{I_{\mathrm{tot}}}\xspace}
\newcommand{\interfVar}{\ensuremath{I}\xspace}

\newcommand{\terminals}{\ensuremath{T}\xspace}
\newcommand{\fixTerminals}{\ensuremath{\terminals_f}\xspace}

\newcommand{\potTerminals}{\ensuremath{\terminals_p}\xspace}

\newcommand{\quota}{\ensuremath{Q}\xspace}
\newcommand{\cost}{\ensuremath{c}\xspace}

\newcommand{\interference}{\ensuremath{I}\xspace}

\newcommand{\setVertices}{\ensuremath{V}\xspace}

\newcommand{\vertex}{\ensuremath{i}\xspace}
\newcommand{\anothervertex}{\ensuremath{j}\xspace}

\newcommand{\vertexCosts}{\ensuremath{b}\xspace}
\newcommand{\vertexProfits}{\ensuremath{q}\xspace}
\newcommand{\vertexCostsNumbered}[1]{\ensuremath{\vertexCosts_{#1}}\xspace}
\newcommand{\vertexProfitsNumbered}[1]{\ensuremath{\vertexProfits_{#1}}\xspace}

\newcommand{\setEdges}{\ensuremath{E}\xspace}
\newcommand{\edge}{\ensuremath{e}\xspace}

\newcommand{\edgeCost}[1]{\ensuremath{\cost_{#1}}\xspace}
\newcommand{\edgeCostVertices}[2]{\ensuremath{\cost_{#1#2}}\xspace}

\newcommand{\graph}{\ensuremath{G}\xspace}
\newcommand{\graphFull}{\ensuremath{\graph = (\setVertices,\setEdges)}\xspace}

\newcommand{\setArcs}{\ensuremath{A}\xspace}
\newcommand{\arc}{\ensuremath{a}\xspace}

\newcommand{\arcCostVertices}[2]{\ensuremath{\cost_{#1#2}}\xspace}
\newcommand{\diGraphFull}{\ensuremath{\ensuremath{D} = (\setVertices,\setArcs)}\xspace}
\newcommand{\incomingArcs}[1]{\ensuremath{\delta^-(#1)}\xspace}
\newcommand{\outgoingArcs}[1]{\ensuremath{\delta^+(#1)}\xspace}

\newcommand{\arcFlow}{\ensuremath{f}\xspace}

\begin{document}

\title{Integrated Wind Farm Design: Optimizing Turbine Placement and Cable Routing with Wake Effects}
\date{}
\author[a]{Jaap~Pedersen\footnote{ORCID: {0000-0003-4047-0042}, corresponding author, pedersen@zib.de}$^,$}
\author[a,b]{Niels~Lindner\footnote{ORCID: {0000-0002-8337-4387}$^,$}}
\author[a,c]{Daniel~Rehfeldt\footnote{ORCID: {0000-0002-2877-074X}}$^,$}
\author[a,d]{Thorsten~Koch\footnote{ORCID: {0000-0002-1967-0077}}$^,$}

\affil[a]{Zuse Institute Berlin, Berlin, Germany}
\affil[b]{Freie Universität Berlin, Berlin, Germany}
\affil[c]{IVU Traffic Technologies, Berlin, Germany}
\affil[d]{Technische Universität Berlin, Berlin, Germany}

\hypersetup{pdftitle={
                Integrated Wind Farm Design: Optimizing Turbine Placement and Cable Routing with Wake Effects},
            pdfauthor={}}

\maketitle
\vspace{-30pt}

\begin{abstract}
An accelerated deployment of renewable energy sources is crucial for a successful transformation of the current energy system, with wind energy playing a key role in this transition. This study addresses the integrated wind farm layout and cable routing problem, a challenging nonlinear optimization problem. We model this problem as an extended version of the quota Steiner tree problem (QSTP), optimizing turbine placement and network connectivity simultaneously to meet specified expansion targets. Our proposed approach accounts for the wake effect $-$ a region of reduced wind speed induced by each installed turbine $-$ and enforces minimum spacing between turbines. We introduce an exact solution framework in terms of the novel quota Steiner tree problem with interference (QSTPI). By leveraging an interference-based splitting strategy, we develop an advanced solver capable of tackling large-scale problem instances. The presented approach outperforms generic state-of-the-art mixed integer programming solvers on our dataset by up to two orders of magnitude. Further, we present a hop-constrained variant of the QSTPI to handle cable capacities in the context of radial topologies. 
Moreover, we demonstrate that our integrated method significantly reduces the costs in contrast to a sequential approach. Thus, we provide a planning tool that enhances existing planning methodologies for supporting a faster and cost-efficient expansion of wind energy.
\end{abstract}

\subsection*{Keywords}
OR in energy, combinatorial optimization, quota Steiner tree problem, wind farm design, wake effect

\section{Introduction}

The transformation of the existing energy system including an increasing deployment of renewable energy sources is vital to reduce global greenhouse gas emissions and to mitigate the impacts of climate change. Wind energy has proven to be an established source of electricity in this transformation \parencite{Veers.2019}. The global wind capacity has increased by more than 20\% annually between 2000 and 2019 \parencite{Pryor.2020}, reaching 730 GW in 2020 \parencite{OurWorldinData.2021}, with a further increase of 50\% expected by the end of 2023 \parencite{Pryor.2020}. Additionally, a sharp decline of the already low cost of wind energy is expected by 2050 \parencite{Jansen.2020, Wiser.2016, Wiser.2021}. This move towards renewable energy sources has created a more competitive market for companies investing in such technologies, making tools for designing and operating such systems optimally vital. With its impressive progress made in the last 20 years \parencite{koch2022}, mathematical programming has proven to be an efficient approach to tackle these challenging problems.

\begin{figure}[ht]
\centering
\begin{subfigure}[t]{0.48\linewidth}
        \centering
        \begin{tikzpicture}[scale=1.5, every node/.style={scale=1.}]
            \tikzstyle{terminal} = [draw, thick, minimum size=0.4, inner sep=2.6pt, color=red];
            \tikzstyle{pot_terminal} = [draw=blue, pattern=north west lines, pattern color=blue, thick, minimum size=0.4, inner sep=2.6pt];
            \tikzstyle{pot_terminald} = [draw, fill, densely dotted, minimum size=0.4, inner sep=2.6pt, color=blue!20];
            \tikzstyle{steiner} = [circle, fill, draw, thick, minimum size=0.2, inner sep=1.5pt];
            \tikzstyle{steinerd} = [circle, densely dotted, draw, minimum size=0.2, inner sep=1.5pt];
            \def\y{16}
            \node[] (yt1) at (0.5+\y, .5) {\includegraphics[width=0.7cm]{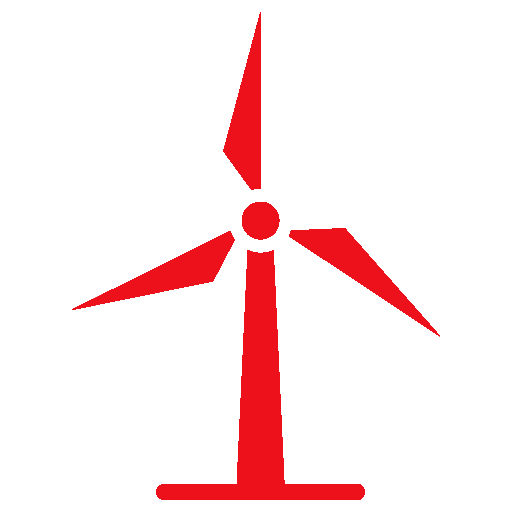}};
            \node[] (yt2) at (0.5+\y, 3.5) {} node[above of=yt2,yshift=-20pt]{\includegraphics[width=0.7cm]{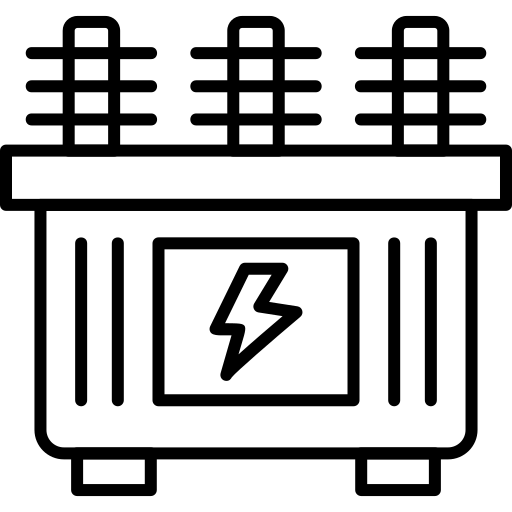}};
            \node[] (ypt1) at(-.5 + \y, 1.) {\includegraphics[width=0.7cm]{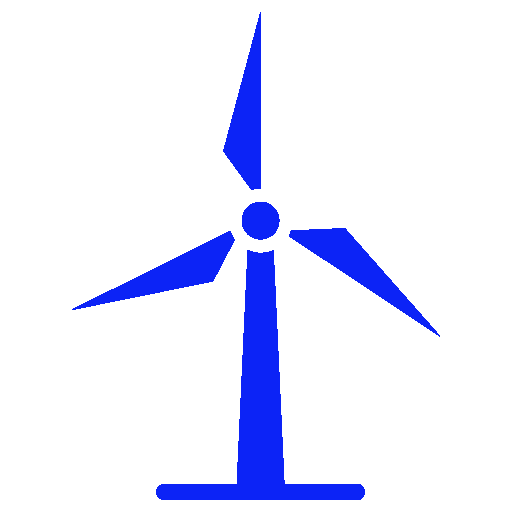}};
            \node[left of=ypt1, xshift=0pt] {\scriptsize 9\,MWh/h};
            \node[] (ypt2) at(1.5 + \y, 1.2) {\includegraphics[width=0.7cm]{turbine.png}};
            \node[right of=ypt2, xshift=0pt] {\scriptsize 8\,MWh/h};
            \node[] (ypt3) at(0. + \y, 2.5) {\includegraphics[width=0.7cm]{turbine.png}};
            \node[left of=ypt3, xshift=0pt] {\scriptsize 5\,MWh/h};
            \node[steiner] (ys1) at (0.5 + \y, 1.3) {};
            \node[steiner] (ys2) at (1.3 + \y, 2.3) {};
                
            \draw[thick, -](yt1) -- node [midway, left=-2.0pt] {} (ys1); 
            
            \draw[thick, -](ys1) -- node [midway, right=-2.0pt] {} (ys2); 
            \draw[thick, -](ys1) -- node [midway, above=-2.0pt] {} (ypt1); 
            \draw[thick, -](ys1) -- node [midway, right=-2.0pt] {} (ypt3); 
            
            \draw[thick, -](ys1) -- node [midway, below=-2.0pt] {} (ypt2); 
            \draw[thick, -](ys2) -- node [midway, above=-2.0pt] {} (ypt3); 
            \draw[thick, -](yt2) -- node [midway, left=-2.0pt] {} (ypt3); 
            
            \draw[thick, -](yt2) -- node [midway, above=-2.0pt] {} (ys2); 
        \end{tikzpicture}
        \caption{IWFLCR instance with substation, three blue potential turbine locations with energy yield, one red fixed turbine, and two black additional nodes for cable routing (icons: \url{flaticon.com}). The aim is to connect the red turbine and some of the black turbines to the substation such that the investment costs are minimized, and a certain energy yield is met.}
    \end{subfigure}
    \hfill
    \begin{subfigure}[t]{0.48\linewidth}
        \centering
        \begin{tikzpicture}[scale=1.5, every node/.style={scale=1.}]
            \tikzstyle{terminal} = [draw, fill, thick, minimum size=0.4, inner sep=2.6pt, color=red];
            \tikzstyle{pot_terminal} = [draw=blue, pattern=north west lines, pattern color=blue, thick, minimum size=0.4, inner sep=2.6pt];
            \tikzstyle{pot_terminald} = [draw, fill, densely dotted, minimum size=0.4, inner sep=2.6pt, color=blue!20];
            \tikzstyle{steiner} = [circle, fill, draw, thick, minimum size=0.2, inner sep=1.5pt];
            \tikzstyle{steinerd} = [circle, densely dotted, draw, minimum size=0.2, inner sep=1.5pt];
            \def\y{16}
            \node[terminal] (yt1) at (0.5+\y, .5) {};
            \node[terminal] (yt2) at (0.5+\y, 3.5) {} node[above of=yt2,yshift=-20pt]{$r$};
            \node[pot_terminal] (ypt1) at(-.5 + \y, 1.) {} node[above of=ypt1, yshift=-20pt]{\scriptsize $i$} node[left of=ypt1, xshift=10pt]{\scriptsize $q_i=9$};
            \node[pot_terminal] (ypt2) at(1.5 + \y, 1.2) {} node[below of=ypt2, yshift=20pt]{\scriptsize $j$} node[right of=ypt2, xshift=-10pt]{\scriptsize $q_j=8$};
            \node[pot_terminal] (ypt3) at(0. + \y, 2.5) {} 
            node[above of=ypt3, xshift=-5pt, yshift=-20pt]{\scriptsize $k$}node[left of=ypt3, xshift=10pt, yshift=-5pt]{\scriptsize $q_k=5$};
            \node[steiner] (ys1) at (0.5 + \y, 1.3) {};
            \node[steiner] (ys2) at (1.3 + \y, 2.3) {};
                
            \draw[thick, -latex](yt1) edge[bend left=20] node [midway, left=-2.0pt] {} (ys1);
            \draw[thick, latex-](yt1) edge[bend left=-20] node [midway, right=-2.0pt] {} (ys1);
            
            \draw[thick, -latex](ys1) edge[bend left=20] node [midway, left=-2.0pt] {} (ys2); 
            \draw[thick, latex-](ys1) edge[bend left=-20] node [midway, right=-2.0pt] {} (ys2);
            \draw[thick, -latex](ys1) edge[bend left=20] node [midway, below=-2.0pt] {} (ypt1);
            \draw[thick, latex-](ys1) edge[bend left=-20] node [midway, above=-2.0pt] {} (ypt1);
            \draw[thick, -latex](ys1) edge[bend left=20] node [midway, left=-2.0pt] {} (ypt3);
            \draw[thick, latex-](ys1) edge[bend left=-20] node [midway, right=-2.0pt] {} (ypt3);
            
            \draw[thick, -latex](ys1) edge[bend left=20] node [midway, above=-2.0pt] {} (ypt2);
            \draw[thick, latex-](ys1) edge[bend left=-20] node [midway, below=-2.0pt] {} (ypt2);
            \draw[thick, -latex](ys2) edge[bend left=20] node [midway, above=-2.0pt] {} (ypt3);
            \draw[thick, latex-](ys2) edge[bend left=-20] node [midway, above=-2.0pt] {} (ypt3);
            \draw[thick, -latex](yt2) edge[bend left=20] node [midway, right=-2.0pt] {} (ypt3);
            \draw[thick, latex-](yt2) edge[bend left=-20] node [midway, left=-2.0pt] {} (ypt3);
            
            \draw[thick, -latex](yt2) edge[bend left=20] node [midway, above=-2.0pt] {} (ys2);
            \draw[thick, latex-](yt2) edge[bend left=-20] node [midway, above=-2.0pt] {} (ys2);  
        \end{tikzpicture}
        \caption{QSTP instance with root $r$, three blue potential terminals $i,j,k$ with their respective quota profits $q_i,q_j,q_k$, one red fixed terminal, and two black Steiner nodes. The aim is to find a minimum cost Steiner arborescence rooted in $r$ and containing the red fixed terminal such that a certain quota is met.}
    \end{subfigure}
    \caption{Translating a wind farm to a combinatorial model.}
    \label{fig:intro}
\end{figure}

In this paper, we investigate the optimal design of wind farms. We present a novel approach for solving the \emph{integrated wind farm layout and cable routing optimization problem} (IWFLCR).
Given a network of potential wind turbine locations and cable routes, the problem is to select turbine positions and a cable routing establishing connections to substations, i.e., the interfaces to the higher-level electricity grid. The goal is to minimize the total investment costs, that may depend on turbine location, turbine type, cable type etc., while fulfilling a certain expansion target in terms of energy yield.
In an abstract language, the IWFLCR can be modeled as a \emph{Quota Steiner tree problem} (\cite{pedersen2024}) as visualized in Figure~\ref{fig:intro}.
This viewpoint however neglects important practical aspects, such as regions of slower wind speed and increased turbulence behind a turbine (wake effects), and bounds on the maximum number of turbines that can be connected by a single cable.

The contributions of our work are as follows. 
For our algorithmic approach, we focus on the combinatorial essence of such a network design problem, and introduce the \emph{Quota Steiner tree problem with interference} (QSTPI). The QSTPI is an extension of the Steiner tree problem in graphs, and is enriched with important technical constraints such as maintaining a minimum distance between turbines and including wake effects. The minimum distances ensures that the blades to not touch and the impact of turbulence is limited \parencite{fischetti2017Phd, fischetti2022}. We therefore offer a discrete optimization model that can both profit from advanced computational techniques and as well as from integrating realistic aspects. We formulate QSTPI as mathematical program in several ways, both as binary quadratic integer program, and as integer linear program. Furthermore, we propose a method to split QSTPI instances in two, exploiting information on lower bounds on the total interference, and thereby accelerating the solution process. Our methods are incorporated into the \scipjack software package \parencite{rehfeldt2021a}, a state-of-the-art branch-and-cut-based solver specialized for Steiner tree problems, moreover advancing its shortest path primal heuristic to handle wake effects. We demonstrate on a realistic benchmark set that our splitting strategy solves more QSTPI instances, and about two orders of magnitude faster, than a black-box approach using the commercial solver \gurobi{} \parencite{gurobi1101}.
Having developed an effective method to solve QSTPI instances, we also provide an extension, called \textit{hop-constrained} QSTPI, that is able to deal with cable capacity limits. To this end, we apply our Steiner tree approach to a layered graph, and enhance our \scipjack implementation with tailored cutting planes. In the end, we can solve hop-constrained QSTPI instances approximately seven times faster than \gurobi{} with a capacitated flow-based model.
Finally, we underline that there is a significant price of sequentiality, i.e., there are noteworthy quantitative benefits that can be obtained by solving the wind farm layout and cable routing optimization simultaneously instead of sequentially. We further analyze the trade-off between using the hop-constrained QSTPI approach on one hand, and using unconstrained QSTPI with an a posteriori capacitated cable routing on the other.

This paper is structured as follows. We begin with a literature review in Section~\ref{sec:literature}.
The QSTPI and its integer programming formulations are presented in Section~\ref{sec:QSTP}. We propose a solution strategy based on splitting the total interference experienced in the wind farm in Section~\ref{sec:Split}. We integrate our proposed model into \scipjack and apply our new methodology in Section~\ref{sec:CompStudy} on a large dataset. In Section~\ref{sec:hop-qstpi}, we extend our formulation to deal with cable capacity constraints. We highlight the benefits of using an integrated approach in contrast to the sequential one in Section~\ref{sec:Res_Sequentiel_vs_combined}. Finally, we discuss and conclude our method and results in Section~\ref{sec:Discussion}.

\section{Literature review}
\label{sec:literature}

We will first review literature on wind farm planning in Section~\ref{sec:literature-wind}, and then proceed in Section~\ref{sec:literature-stp} with our mathematical background, the Steiner tree problem.

\subsection{Literature on wind farm planning}
\label{sec:literature-wind}

Designing wind farms involves a variety of decisions leading to many challenging optimization problems. In this paper, we focus on the integration of two main design tasks: the wind farm layout optimization (WFLO) and the wind farm cable routing (WFCR). There exists a vast amount of literature investigating both the WFLO and the WFCR. For an extensive overview we refer the reader to the surveys by \textcite{Samorani2013}, \textcite{hou2019review}, and \textcite{fischetti2019overview}.

The WFLO maximizes the energy output of the wind farm. It does not only consider the number of turbines to be installed, but also takes the \textit{wake effect} into account. A region of slower wind speed and increasing turbulence, a \textit{wake}, is caused by each installed turbine, resulting in a reduced energy yield of turbines downwind. An investigation of the offshore wind farm areas along the US east coast by \textcite{pryor2021wind} shows that power reduction can be reduced by one-third due to wakes caused by upwind turbines and wind farms. Therefore, it is inevitable to consider the wake effect when modeling the WFLO. There exists a variety of studies modeling the WFLO as a mixed integer programming problem (MIP). A first MIP is presented in 
\textcite{donovan2005wind} based on the independent set problem. \textcite{turner2014new} approximate the nonlinear optimization problem of minimizing the wake effect by a quadratic integer and a mixed integer linear programming model (MILP). Based on the MILP proposed by \textcite{archer2011}, \textcite{fischetti2016proximity} combine the MIP approach with a proximity search heuristic to solve large-scale instances with the drawback of losing proven global optimality. \textcite{fleming2016wind} and \textcite{gebraad2017maximization} couple the positioning and controlling of the turbines in the planning phase to furhter increase the annual energy production.

In the WFCR, the costs of connecting the chosen turbines to the power grid are minimized. In the case of offshore wind farms, \textcite{gonzalez2014} point out that the costs of the electrical infrastructure account for 15\% to 30\% of the overall initial costs of the wind farm, making them an essential matter of expenses. During the operation, power losses play a crucial role in terms of expected revenue. For example, \textcite{Fischetti.2018} present a combination of MILP formulation and metaheuristics, named \textit{matheuristics}, to minimize not only the initial costs, but also the reduced revenue due to power losses in the lifetime of a wind farm by considering cable capacities.

Due to the complexity of the WFLO and WFCR, the problems are usually solved sequentially. However, by design these problems are conflicting, i.e., WFLO aims to spread the turbines as far as possible to reduce wakes, whereas the WFCR tries to keep them as close as possible to reduce cable costs. Recent promising studies have been conducted to integrate these two problems balancing these conflicts \parencite{Fischetti.2021a, fischetti2022, cazzaro2023}. \textcite{fischetti2022} model the integrated layout and cable routing problem as a mixed-integer linear program in the context of offshore wind farms. The authors introduce a novel set of cutting planes and present combinatorial separation procedures integrated in a branch-and-cut solver. Although the problem is solved exactly, the approach is only performed on small instances with up to 40 potential turbine positions, whereas in practice, large wind farms might include up to 100~--~200 turbines\footnote{\url{https://en.wikipedia.org/wiki/List_of_offshore_wind_farms}; Retrieved on 01/08/2025}. Recently, \textcite{cazzaro2023} approach the integrated layout and cable optimization problem heuristically using an adapted variable neighborhood search. Even though the heuristic combining layout and cable routing provides better solutions than the previously proposed sequential approach, it cannot provide any conclusion in terms of quality with respect to a globally optimal solution.

\subsection{Literature, algorithms and software for Steiner trees}
\label{sec:literature-stp}

The goal of this study is to provide a novel approach for solving the IWFLCR considering wake effects, minimum distance between turbines, and balancing investment cost and revenues. Conceptually, the IWFLCR is a network design problem with a number of additional constraints giving rise to a tree-shaped topology \parencite{Fischetti.2018}: We model the IWFLCR as an extension to the Steiner tree problem in graphs (STP), a classical NP-hard combinatorial problem \parencite{karp1972}, and one of the most studied problems in combinatorial optimization. Given an undirected graph with non-negative edge costs, the STP aims to find a tree that interconnects a given set of special points, referred to as \emph{terminals}, at minimum cost. The STP and its variations arise in many real-world applications like, e.g., network design problems in telecommunication, electricity, or in district heating, as well as other fields such as biology; see, e.g., \textcite{leitner2014}, \textcite{bolukbasi2018}, \textcite{ljubic2006}, and \textcite{klimm2020}, respectively. For a comprehensive overview of the STP and its variant the readers are referred to a recent surveys by \textcite{ljubic2021}, \textcite{rehfeldt2021a}, and \textcite{Pedersen2024Prize}, and the references therein. In the PhD thesis by \textcite{ridremont2019}, the STP is applied to find a robust cable network for the purpose of designing a wind farm. Recently, \textcite{pedersen2024} have investigated the quota-constrained variant of the STP (QSTP) in the context of large-scale onshore wind farm planning. However, in the context of Steiner trees, vital technical details, such as, among others, the wake effect or the minimum distance between turbines, are only implicitly modeled in the input data rather than in the optimization model. The authors introduce a novel transformation of the QSTP, which is integrated into the specialized Steiner tree solver \scipjack, vastly outperforming state-of-the-art general out-of-the-box solvers. As an extension of the non-commercial general MIP solver \scip \parencite{bestuzheva2023}, \scipjack uses a branch-and-cut procedure to handle the exponential number of constraints induced by the Steiner cut-like constraints. These constraints will be recalled in Section~\ref{sec:DCutForm}. By additionally using various problem-specific reduction techniques and heuristics, \scipjack has proven that it efficiently solves STP-related problems of magnitudes larger than what general out-of-the-box MIP solvers can even load into memory; see, e.g., \textcite{rehfeldt2019, rehfeldt2022, RehfeldtKoch2023}. For a comprehensive description of \scipjack's operational principle the reader is referred to the PhD thesis by~\textcite{rehfeldt2021a}.

\section{The Quota Steiner tree problem}\label{sec:QSTP}

We model the integrated wind farm layout and cable routing problem as a \emph{Quota Steiner tree problem with interference} (QSTPI), extending and reformulating the classical quota Steiner tree problem (\cite{johnson2000}). In particular, we consider arbitrary tree-shaped network topologies. We introduce the QSTPI formally in Section~\ref{sec:qstpi}, present a preliminary integer programming (IP) formulation in terms of directed cuts in Section~\ref{sec:DCutForm}, before  applying the techniques by \cite{pedersen2024} to obtain our final IP model in Section~\ref{sec:qstpTrans}.

\subsection{The Quota Steiner tree problem with interference}
\label{sec:qstpi}

We adapt the formulation by \textcite{johnson2000} of the prize-collecting quota Steiner tree problem by adding vertex costs and integrating interference effects into the profit.
Let \graphFull be an undirected graph, whose vertex set $\setVertices$ is partitioned into a set of \textit{fixed terminals} $\fixTerminals$, a set of \textit{potential terminals} $\potTerminals$, and a set of additional \textit{Steiner nodes}. We associate the potential terminals with \emph{vertex costs} $\vertexCosts : \potTerminals\rightarrow \R_{>0}$ and \emph{quota profits} $\vertexProfits: \potTerminals \rightarrow \R_{>0}$, referring to a \emph{quota} $Q \in \R_{>0}$. Moreover, \emph{edge costs} $\cost : \setEdges \rightarrow \R_{\ge0}$ are given.

In terms of wind farm planning, the substations of the grid constitute $\fixTerminals$, while $\potTerminals$ represents the available wind turbine positions. The Steiner nodes allow for a more flexible cable routing. The vertex costs $\vertexCosts$ represent the investment costs for a turbine, e.g., foundation costs, which depend on the potential position, and the cost of the turbine itself, which may vary based on the chosen turbine type. The quota profits $\vertexProfits$ the turbine's annual energy yield under a certain wind scenario. The built turbines should contribute an output of at least the quota $Q$.
The set of edges $\setEdges$ represents the cable connections we can choose from, and $\cost$ comprises their costs.

In extension to the QSTP presented by \textcite{pedersen2024}, each potential terminal, if chosen, may interfere with other potential terminals, reducing their assigned quota profit. As to wind turbine placement, each built turbine introduces a \textit{wake}, an area of slower wind speed downstream, decreasing the power production of the following turbines, and, consequently, reducing its annually energy yield. We pin down this \textit{wake effect} in a twofold way by \emph{interference} and \emph{minimum distance} notions: Let $\interference_{ij}\in \R_{\ge0}$ be the interference experienced by the potential terminal ${j} \in \potTerminals$ if the potential terminal ${i} \in \potTerminals$ is chosen. Furthermore, let $d_{ij}\in \R_{\geq 0}$ denote the \emph{distance} between two potential terminals $i, j \in T_p$, and let $D_{\mathrm{min}}$ be a \emph{minimum distance}. While interference and minimum distance have been discussed before by \textcite{fischetti2016proximity, fischetti2022}, we integrate these concepts into the QSTP from \textcite{pedersen2024}.

\begin{definition}[Quota Steiner tree problem with Interference]
Given an instance \linebreak $(G, T_f, T_p, b, q, c, I, Q, d, D_{\mathrm{min}})$ with $G = (V, E)$ as above, the \emph{Quota Steiner tree problem with Interference (QSTPI)} is to find a tree $S = (\setVertices^\prime, \setEdges ^\prime) \subseteq \graph$ that contains all fixed terminals $\fixTerminals$, minimizes the total cable and turbine costs
\begin{equation}
    C(S) \coloneqq \sum_{(i,j)\in \setEdges^\prime} \edgeCostVertices{ij}{} + \sum_{i \in \potTerminals \cap \setVertices^\prime} \vertexCostsNumbered{i},
\end{equation}
fulfills the quota taking the induced interference into account:
\begin{equation}
    \quota(S) \coloneqq \sum_{i \in \potTerminals \cap \setVertices^\prime}\left( \vertexProfitsNumbered{i} -  \sum_{j \in \potTerminals \cap \setVertices^\prime \setminus \{i\}} \interference_{ji}\right) \ge \quota,\label{eq:quota_definition}
\end{equation}
and guarantees that $d_{ij} \geq D_{\mathrm{min}}$ for all $i, j \in \setVertices^\prime \cap \potTerminals$ with $i \neq j$.
\end{definition}

The QSTPI can also be seen as a single-objective specialization of a multi-objective problem, where cable and turbine costs are supposed to be minimized while simultaneously maximizing the output. In the QSTPI, a minimum output level is ensured by \eqref{eq:quota_definition} and the number of turbines depends on the chosen quota, their individual quota profit, and the total induced interference, and, thus, is not known beforehand.

\subsection{Directed cut formulation}\label{sec:DCutForm}

To develop an integer programming (IP) formulation, we transform the QSTPI into a variant of the \textit{Steiner arborescence problem} (SAP). The general SAP is presented, e.g., in \textcite{wong1984}. The original undirected graph $\graph$ is transformed into a directed graph $\diGraphFull$ where $\setArcs\coloneqq \lbrace (i,j),(j,i):\{i,j\} \in \setEdges\rbrace$. By shifting the costs of a vertex $\vertex$ onto the costs of its incoming arcs (see \cite{ljubic2006}), the arc costs are defined as
\begin{align}
& \arcCostVertices{i}{j} \coloneqq \begin{dcases} 
        \edgeCost{e} + \vertexCostsNumbered{j} & \text{if } j \in \potTerminals,\\
        \edgeCost{e}& \text{otherwise},
        \end{dcases} & \forall \arc=(i,j) \in \setArcs, \label{eq:cost_shift}
\end{align}
where $\edgeCost{e}$ represents the cost of the corresponding edge $\edge=\{i,j\}$ in the original undirected graph $\graph$. For a subset of vertices $W\subseteq V$, we denote $\delta^+(W) \coloneqq \lbrace (i,j) \in \setArcs: i \in W, j \in \setVertices\backslash W \rbrace$ as the set of outgoing arcs and $\delta^-(W) \coloneqq \lbrace (i,j) \in A: i \in \setVertices\backslash W, j \in W \rbrace$ as the set of incoming arcs. For a single vertex ${i}$, we write $\delta^+({i}) \coloneqq \delta^+(\lbrace {i}\rbrace)$ and $\delta^-({i}) \coloneqq \delta^-(\lbrace {i}\rbrace)$. For any set $K$ and any function $x:K \to \R$, we define $x(K) \coloneqq \sum_{i \in K} x_i$. 

We introduce a binary variable $\arcVar_{ij}$ for each $(i,j) \in \setArcs$ if the arc $(i,j)$ is contained in the Steiner tree ($\arcVar_{ij} = 1$) or not ($\arcVar_{ij} = 0$). Furthermore, let $\nodeVar_k$ be a binary variable for each $k \in \potTerminals$ indicating whether or not the potential terminal $k$ is chosen. The rooted directed cut IP formulation of the QSTPI with an arbitrary root $r \in \fixTerminals$ is given as follows:
\begin{align}
    \min \quad&\cost^T \arcVar &\label{eq:obj_qstp}\\
    \text{s.t.} \quad
    & \arcVar(\delta^-(W)) \ge 1& \forall W \subset \setVertices, r\notin W, |W \cap \fixTerminals| \ge 1\label{Steiner cut}\\
    & \arcVar(\delta^-(W)) \ge \nodeVar_i& \forall W \subset \setVertices, r\notin W, |W \cap \potTerminals| \ge 1,  {i}\in \potTerminals \label{Steiner cut pot}\\
    & \sum_{{i} \in \potTerminals} \left(\vertexProfitsNumbered{i} - \sum_{{j}\in \potTerminals} \interference_{ji} \nodeVar_j\right) \nodeVar_i\ge \quota & \label{quota cons}\\
    & \nodeVar_{i} + \nodeVar_{j} \le 1 & \forall i, j \in \potTerminals, i \neq j, d_{ij} < D_{\mathrm{min}} \label{eq:minD_qstp}\\
    & \arcVar_{ij}, \nodeVar_{k}  \in \lbrace 0,1 \rbrace &\forall (i,j) \in \setArcs, \forall {k} \in \potTerminals \label{eq:var_qstp}
\end{align}
The Steiner cut constraint \eqref{Steiner cut} requires that any subset of nodes of the graph $W \subset \setVertices$ that does not contain the root ($r\notin W$), but at least one fixed terminal, i.e., $|W \cap \fixTerminals| \ge 1$, has at least one incoming arc. In other words, for any cut between the root node and a fixed terminal, there has to exists a least one arc crossing this cut towards the side of the fixed terminal. This guarantees that there exists a path from the root $r$ to each $t\in \fixTerminals$. Note, that if there exists multiple substations, i.e, $|\fixTerminals|>1$, then there would exist a connection between the substations which is not desirable. However, by inserting an artificial root node $r^\prime$ connected to each substation $t\in\fixTerminals$ with an artificial arc $(r^\prime, t)$ with zero costs, one makes sure that there exists no directed path between any two substations. The artificial root node can also be seen as the higher level grid. \eqref{Steiner cut pot} guarantees that there exists a path from $r$ to the potential terminal $i\in \potTerminals$ if it is chosen to contribute to the quota constraint \eqref{quota cons} ($\nodeVar_i = 1$). Finally, \eqref{eq:minD_qstp} ensures that no two potential terminals are chosen that are too close to each other.

\subsection{Transformation}
\label{sec:qstpTrans}

Having modeled the QSTPI as an IP by means of the directed cut formulation, we now shortly revisit the transformation of QSTP presented in \textcite{pedersen2024}, which we extend by interference and minimal distance constraints. For each potential terminal $i \in \potTerminals$, we add a new fixed terminal $i^\prime \in \fixTerminals$. For each of these new fixed terminals ${i^\prime}$, we create an arc $(r,{i^\prime})$ with costs $\arcCostVertices{r}{{i^\prime}} \coloneqq 0$ and an arc $(i, {i^\prime})$ with costs $\arcCostVertices{i}{{i^\prime}} \coloneqq 0$, and add them to $\setArcs$; see Fig.~\ref{fig:GraphTrans}.
Since $i^\prime$ is a fixed terminal, the Steiner cut constraint \eqref{eq:Steiner cut trans} then requires that $\arcVar_{ii^\prime} = 1$ or $\arcVar_{ri^\prime} = 1$. We interpret the first case as turbine $i \in \potTerminals$ being built, and the latter case as turbine $i$ not being built.
Following the argumentation in \textcite{pedersen2024}, instead of counting the collected quota, we limit the amount of quota we can \textit{waste}.
Quota is wasted when a turbine $i$ is not built, so that $\arcVar_{ri^\prime} = 1$, see \eqref{eq:quota_cons_trans}.
The total interference $\interfTotVar$ is given by the interference between all chosen turbines,
see \eqref{eq:interf_tot_trans}. The minimum distance constraint \eqref{eq:minD_trans} prevents building simultaneously two turbines that are too close. Note, that although the Steiner cut constraint~\eqref{eq:Steiner cut trans} allows for solutions with $x(\incomingArcs{i^\prime}) = \arcVar_{ii^\prime} + \arcVar_{ri^\prime} = 2$, setting $\arcVar_{ii^\prime} \coloneqq 0$ and $\arcVar_{ri^\prime} \coloneqq 1$ results in a feasible solution with the same or less cost.

The transformed QSTP including interference and minimum distance constraint is given as follows:
\begin{align}
    \text{QSTPI-TRANS}: &&\nonumber\\
    \min \quad&\cost^T \arcVar &\label{eq:obj_trans_qstp}\\
    \text{s.t.} \quad
    & \arcVar(\delta^-(W)) \ge 1& \forall W \subset \setVertices, r\notin W, |W \cap \fixTerminals| \ge 1\label{eq:Steiner cut trans}\\
    & \sum_{i \in \potTerminals} \vertexProfitsNumbered{i} \arcVar_{r i^{\prime}} + \interfTotVar \le \sum_{i \in \potTerminals} \vertexProfitsNumbered{i} - \quota & \label{eq:quota_cons_trans}\\
    &\interfTotVar \ge \sum_{i \in \potTerminals} \sum_{j \neq i \in \potTerminals} \interference_{ij} \arcVar_{i i^{\prime}} \arcVar_{j j^{\prime}}&\label{eq:interf_tot_trans}\\
    & \arcVar_{i i^{\prime}} + \arcVar_{j j^{\prime}} \le 1 & \forall i, j \in \potTerminals, i \neq j, d_{ij} < D_{\mathrm{min}}\label{eq:minD_trans}\\
    & \arcVar_{i,j} \in \lbrace 0,1 \rbrace &\forall (i,j) \in \setArcs\label{eq:var_trans_qstp}
\end{align}

The formulation QSTPI-TRANS is a quadratic-constrained binary optimization problem due to \eqref{eq:interf_tot_trans} and NP-hard \parencite{pedersen2024}. Due to the exponential number of  constraints induced by the Steiner cut constraints \eqref{eq:quota_cons_trans}, it is reasonable to use a solver which is specialized in separating these constraints, e.g., \scipjack\parencite{RehfeldtKoch2023}. As an extension to the general framework \scip\parencite{bestuzheva2023}, which can handle mixed integer nonlinear programs (MINLP), one could integrate the quadratic constraint \eqref{eq:interf_tot_trans} directly. However, preliminary experiments have shown that explicitly linearizing the bilinear products $\arcVar_{ii^\prime}\arcVar_{jj^\prime}$ in \eqref{eq:interf_tot_trans} using McCormick's linearization \parencite{mccormick1976computability} appears to be faster than let \scip handle the nonlinearities. By introducing additionally variables $\bilinVar_{ij}$, \eqref{eq:interf_tot_trans} can be replaced as follows:
\begin{align}
    &\interfTotVar \ge \sum_{i \in \potTerminals}\sum_{j \neq i \in \potTerminals} \interference_{ij} \bilinVar_{ij}&\label{eq:quota_cons_trans_lin}\\
    &\bilinVar_{ij} \ge \arcVar_{i i^{\prime}} + \arcVar_{j j^{\prime}} - 1& \forall i < j \in \potTerminals,\label{eq:bilinear_lin}\\
    &\bilinVar_{ij} \ge 0 & \forall i < j \in \potTerminals.\label{eq:bilinear_var}
\end{align}
As the total interference $\interfTotVar$ is bounded from above in \eqref{eq:quota_cons_trans}, no upper bound has to be imposed on $\bilinVar$. 

\begin{figure}
    \begin{subfigure}{0.4\linewidth}
        \centering
        \begin{tikzpicture}[scale=1.5, every node/.style={scale=1.}]
            \tikzstyle{terminal} = [draw, fill, thick, minimum size=0.4, inner sep=2.6pt, color=red];
            \tikzstyle{pot_terminal} = [draw=blue, pattern=north west lines, pattern color=blue, thick, minimum size=0.4, inner sep=2.6pt];
            \tikzstyle{pot_terminald} = [draw, fill, densely dotted, minimum size=0.4, inner sep=2.6pt, color=blue!20];
            \tikzstyle{steiner} = [circle, fill, draw, thick, minimum size=0.2, inner sep=1.5pt];
            \tikzstyle{steinerd} = [circle, densely dotted, draw, minimum size=0.2, inner sep=1.5pt];
            \def\y{16}
            \node[terminal] (yt1) at (0.5+\y, .5) {};
            \node[terminal] (yt2) at (0.5+\y, 3.5) {} node[above of=yt2,yshift=-20pt]{$r$};
            \node[pot_terminal] (ypt1) at(-.5 + \y, 1.) {} node[above of=ypt1, yshift=-20pt]{\scriptsize $i$} node[left of=ypt1, xshift=10pt]{\scriptsize $q_i=9$};
            \node[pot_terminal] (ypt2) at(1.5 + \y, 1.2) {} node[below of=ypt2, yshift=20pt]{\scriptsize $j$} node[right of=ypt2, xshift=-10pt]{\scriptsize $q_j=8$};
            \node[pot_terminal] (ypt3) at(0. + \y, 2.5) {} 
            node[above of=ypt3, xshift=-5pt, yshift=-20pt]{\scriptsize $k$}node[left of=ypt3, xshift=10pt, yshift=-5pt]{\scriptsize $q_k=5$};
            \node[steiner] (ys1) at (0.5 + \y, 1.3) {};
            \node[steiner] (ys2) at (1.3 + \y, 2.3) {};
                
            \draw[thick, -latex](yt1) edge[bend left=20] node [midway, left=-2.0pt] {} (ys1);
            \draw[thick, latex-](yt1) edge[bend left=-20] node [midway, right=-2.0pt] {} (ys1);
            
            \draw[thick, -latex](ys1) edge[bend left=20] node [midway, left=-2.0pt] {} (ys2);
            \draw[thick, latex-](ys1) edge[bend left=-20] node [midway, right=-2.0pt] {} (ys2);
            \draw[thick, -latex](ys1) edge[bend left=20] node [midway, below=-2.0pt] {} (ypt1);
            \draw[thick, latex-](ys1) edge[bend left=-20] node [midway, above=-2.0pt] {} (ypt1);
            \draw[thick, -latex](ys1) edge[bend left=20] node [midway, left=-2.0pt] {} (ypt3);
            \draw[thick, latex-](ys1) edge[bend left=-20] node [midway, right=-2.0pt] {} (ypt3);
            
            \draw[thick, -latex](ys1) edge[bend left=20] node [midway, above=-2.0pt] {} (ypt2);
            \draw[thick, latex-](ys1) edge[bend left=-20] node [midway, below=-2.0pt] {} (ypt2);
            \draw[thick, -latex](ys2) edge[bend left=20] node [midway, above=-2.0pt] {} (ypt3);
            \draw[thick, latex-](ys2) edge[bend left=-20] node [midway, above=-2.0pt] {} (ypt3);
            \draw[thick, -latex](yt2) edge[bend left=20] node [midway, right=-2.0pt] {} (ypt3);
            \draw[thick, latex-](yt2) edge[bend left=-20] node [midway, left=-2.0pt] {} (ypt3);
            
            \draw[thick, -latex](yt2) edge[bend left=20] node [midway, above=-2.0pt] {} (ys2);
            \draw[thick, latex-](yt2) edge[bend left=-20] node [midway, above=-2.0pt] {} (ys2);  
        \end{tikzpicture}
        \caption{Original instance.}
    \end{subfigure}
    \begin{subfigure}{0.58\linewidth}
    \centering
        \begin{tikzpicture}[scale=1.5, every node/.style={scale=1.}]
            \tikzstyle{terminal} = [draw, fill, thick, minimum size=0.4, inner sep=2.6pt, color=red];
            \tikzstyle{pot_terminal} = [draw=blue, pattern=north west lines, pattern color=blue,
            thick, minimum size=0.4, inner sep=2.6pt];
            \tikzstyle{pot_terminald} = [draw, fill, densely dotted, minimum size=0.4, inner sep=2.6pt, color=blue!20];
            \tikzstyle{steiner} = [circle, fill, draw, thick, minimum size=0.2, inner sep=1.5pt];
            \tikzstyle{steinerd} = [circle, densely dotted, draw, minimum size=0.2, inner sep=1.5pt];
            \def\y{16}
            \node[terminal] (yt1) at (0.5+\y, .5) {};
            \node[terminal] (yt2) at (0.5+\y, 3.5) {} node[above of=yt2,yshift=-20pt]{$r$};
            \node[pot_terminal] (ypt1) at(-.5 + \y, 1.) {} node[above of=ypt1, yshift=-20pt]{\scriptsize $i$};
            \node[pot_terminal] (ypt2) at(1.5 + \y, 1.2) {} node[below of=ypt2, yshift=20pt]{\scriptsize $j$};
            \node[pot_terminal] (ypt3) at(0. + \y, 2.5) {} node[above of=ypt3, xshift=-5pt, yshift=-20pt]{\scriptsize $k$};
            \node[steiner] (ys1) at (0.5 + \y, 1.3) {};
            \node[steiner] (ys2) at (1.3 + \y, 2.3) {};
            \draw[thick, -latex](yt1) edge[bend left=20]   (ys1); 
            \draw[thick, latex-](yt1) edge[bend left=-20]  (ys1); 
            \draw[thick, -latex](ys1) edge[bend left=20]   (ys2); 
            \draw[thick, latex-](ys1) edge[bend left=-20]  (ys2); 
            \draw[thick, -latex](ys1) edge[bend left=20]   (ypt1);
            \draw[thick, latex-](ys1) edge[bend left=-20]  (ypt1);
            \draw[thick, -latex](ys1) edge[bend left=20]   (ypt3);
            \draw[thick, latex-](ys1) edge[bend left=-20]  (ypt3);
            \draw[thick, -latex](ys1) edge[bend left=20]   (ypt2);
            \draw[thick, latex-](ys1) edge[bend left=-20]  (ypt2);
            \draw[thick, -latex](ys2) edge[bend left=20]   (ypt3);
            \draw[thick, latex-](ys2) edge[bend left=-20]  (ypt3);
            \draw[thick, -latex](yt2) edge[bend left=20]   (ypt3);
            \draw[thick, latex-](yt2) edge[bend left=-20]  (ypt3);
            \draw[thick, -latex](yt2) edge[bend left=20]   (ys2); 
            \draw[thick, latex-](yt2) edge[bend left=-20]  (ys2); 
            \draw[thick, latex-](yt2) edge[bend left=-20]  (ys2); 
            \node[terminal](yt11) [left of=ypt1, xshift=-5pt] {} node[above of=yt11, xshift=5pt, yshift=-20pt]{\scriptsize$i'$};
            \node[terminal](yt22) [right of=ypt2, xshift=5pt] {} node[below of=yt22, yshift=20pt]{\scriptsize$j'$};
            \node[terminal](yt33) [left of=ypt3, xshift=-5pt] {} node[above of=yt33, yshift=-20pt]{\scriptsize$k'$};
            \draw[thick, -latex](yt2) edge[bend left=-60] node [midway, left=2pt]{\scriptsize \textbf{0~}}  (yt11);
            \draw[thick, -latex](ypt1) edge node [midway, below=2pt]{\scriptsize \textbf{0~}} (yt11);
            \draw[thick, -latex](yt2) edge[bend left=60] node [midway, right=2pt]{\scriptsize \textbf{0~}}  (yt22);
            \draw[thick, -latex](ypt2) edge node [midway, below=2pt]{\scriptsize \textbf{0~}} (yt22);
            \draw[thick, -latex](yt2) edge[bend left=-5] node [midway, left=2pt]{\scriptsize \textbf{0~}}  (yt33);
            \draw[thick, -latex](ypt3) edge node [midway, below=2pt]{\scriptsize \textbf{0~}} (yt33);
        \end{tikzpicture}
        \caption{Transformed instance.}
    \end{subfigure}
    \caption{Transformation of QSTPI. For each potential terminal $i \in \potTerminals$ (blue squares) a fixed terminal (red square) $i^\prime$ is added, which is connected to the root node $r$ and the original vertex by an arc of zero costs; Steiner nodes are shown as filled circles.} 
    \label{fig:GraphTrans}
\end{figure}

\subsection{The flow-based approach}\label{sec:flow}

Due to the exponential number of constraints \eqref{eq:Steiner cut trans}, which cannot be handled by general MIP solvers out-of-the-box, the problem can alternatively be modeled by a standard flow-based MIP formulation, commonly used in such applications; see, e.g., \textcite{fischetti2022}. The flow-based MIP formulation of the non-transformed QSTPI is given by
\begin{align}
\mathrm{FLOW:}&&\nonumber\\
        \min \quad&\cost^T \arcVar& \label{eq:obj_flow}\\
    \text{s.t.} \quad
    &\sum_{\vertex\in\potTerminals} \vertexProfitsNumbered{\vertex} \nodeVar_\vertex  - \interfTotVar\ge \quota \label{eq:quota_flow}\\
    &\interfTotVar \ge \sum_{\vertex \in \potTerminals}\sum_{\anothervertex \neq \vertex \in \potTerminals} \interference_{\anothervertex\vertex} \nodeVar_\anothervertex \nodeVar_\vertex\label{eq:interference_flow}\\
    &\sum_{\arc\in\incomingArcs{\vertex}}\arcFlow_\arc - \sum_{\arc\in\outgoingArcs{\vertex}}\arcFlow_\arc = \begin{dcases} 0 \quad \forall \vertex \in \setVertices\setminus(\fixTerminals \cup\potTerminals)\\ 
    1 \quad \forall \vertex \in \fixTerminals\setminus \{r\}\\
    \nodeVar_\vertex \quad \forall \vertex \in \potTerminals
    \end{dcases} \label{eq:flowBalance_flow}\\
    & \nodeVar_{i} + \nodeVar_{j} \le 1 & \forall i \neq j \in \potTerminals, d_{ij} < D_{\mathrm{min}}\label{eq:minD_flow}\\
    &\arcVariable_\arc \le y_\vertex &\forall \arc \in \incomingArcs{\vertex}, \forall \vertex \in \potTerminals \label{eq:activeTerm_flow}\\
    &\arcFlow_\arc \le M \arcVar_\arc &\forall \arc \in \setArcs \label{eq:activeArc_flow}\\
    &\nodeVar_\vertex \in \lbrace 0,1\rbrace &\forall\vertex \in \potTerminals\\
    &\arcVar_\arc \in \lbrace 0,1\rbrace, \,\arcFlow_\arc \in \R_{\ge0} &\forall\arc \in \setArcs \label{eq:var_flow}
\end{align}
where $\arcVar$ and $\nodeVar$ denote the decision variables if an arc and a vertex is chosen, respectively, and $\arcFlow$ describes the flow over the arcs. The costs are minimized in \eqref{eq:obj_flow}, note, that the costs $c$ contain the turbine costs if chosen, see \eqref{eq:cost_shift}. Constraint \eqref{eq:quota_flow} describes the quota constraint including the total inference $\interfTotVar$ caused by the chosen turbines in \eqref{eq:interference_flow}. The flow balance at each vertex depends on its type and is captured by constraint \eqref{eq:flowBalance_flow}.  For any Steiner node $\vertex \in \setVertices\setminus(\fixTerminals\cup\potTerminals)$ the ingoing and outgoing flow has to be equal. Any fixed terminal, i.e., the substations, which is not the root, has to be connected to the root node, so a ``flow'' $\arcFlow$ of one has to reach the terminal, i.e., the right-hand side is equal to one. Also, any chosen potential terminal (turbine), i.e., $\nodeVar_\vertex = 1$, has to be connected to the root node. The incoming arcs of a potential terminal can only be active if the potential terminal is chosen, as in \eqref{eq:activeTerm_flow}. Equation \eqref{eq:activeArc_flow} ensures that a flow over an arc is only possible if the arc is active. The big-M notation is used to limit the flow on an active arc, i.e., the upper bound on $\arcFlow_\arc$. We choose $M = | \fixTerminals \cup \potTerminals |$ as an upper limit, which would allow all fixed and potential terminals to be connected via a single string.

\section{Split the problem in two}
\label{sec:Split}

In this section, we introduce a splitting strategy based on lower bounds on the interference. We will use this method to accelerate the solution process for practical instances of the QSTPI later on.

\subsection{Bounding local interference}
Concerning wind farms, the layout problem promotes solutions with wide-spread wind turbine positions, whereas the cable routing problem aims to have them as close as possible to the substation. Furthermore, the linearization \eqref{eq:bilinear_lin} allows for highly fractional solutions of the LP-relaxation of the QSTPI-TRANS (\eqref{eq:obj_trans_qstp}--\eqref{eq:quota_cons_trans} and \eqref{eq:var_trans_qstp} -- \eqref{eq:bilinear_var}), i.e., relaxation \eqref{eq:var_trans_qstp} to $x_{ij}\in [0, 1] \forall(i,j) \in A$, with values $x_{ii^\prime} = x_{jj^\prime} = 0.5$ and $m_{ij} = 0$, which means that $i$ and $j$ contribute to the quota, but not to the interference value, thus causing highly fractional solutions. \textcite{fischetti2022} present cuts for the flow-based formulation that increase the lower bound of interference caused by choosing to install a turbine, and that have also been useful for some hard quadratic problems, including the Quadratic Assignment Problem, see \textcite{fischetti2012three}. A particular and effective cut is given by
\begin{align}
    \sum_{j\neq i \in \potTerminals}\interfVar_{ji} \arcVar_{jj^\prime} \arcVar_{ii^\prime}\ge \mathrm{LB}_i \arcVar_{ii^\prime},&\quad \forall i \in \potTerminals,
\end{align}
i.e., the total interference caused by turbine $i$ if built is bounded by a value $\mathrm{LB}_i$.
The authors propose that if one aims to build $N_{\mathrm{min}}$ turbines, then $\mathrm{LB}_i$ can be chosen as the sum of the $N_{\mathrm{min}} - 1$ lowest values of $\interference_{ij}$ for all $i,j$ with $d_{ij}\ge D_{\mathrm{min}}$. One can extended this thought towards the total induced interference by taking the sum of the $N_{\mathrm{min}}$ lowest values of $\mathrm{LB}_i$ to determine a lower bound $\interfVar_{\mathrm{LB}}$: Let  $(\mathrm{LB}_{(1)}, \dots, \mathrm{LB}_{(|\potTerminals|)})$ be the ascending sort of all $\mathrm{LB}_i$, i.e. $\mathrm{LB}_{(1)}\le\dots \leq \mathrm{LB}_{(|\potTerminals|)}$. Then the lower bound of the total induced interference is given by:
\begin{align}
    \interfTotVar \ge I_{\mathrm{LB}} = \sum_1^{N_{\mathrm{min}}} \mathrm{LB}_{(i)}.
\end{align}

\subsection{Splitting the total interference}

Considering large areas with a high number of potential wind turbine location from which only a small subset is chosen to be built, the previous described lower bound easily is small compared to the actual induced interference. Motivated by this fact, we split QSTPI-TRANS into two subproblems as follows: We define a split value $I_{\mathrm{split}}\in \R_{>0}$ for the total interference, which will be determined in Section~\ref{sec:MinI_lowerbound}. In the first subproblem $P_{\ge}$, which is called the \textit{upproblem}, this split value is imposed as a lower bound on the total interference, i.e., $\interfTotVar \ge I_{\mathrm{split}}$. In the second subproblem $P_{\le}$, the \textit{downproblem}, the split value represents the upper bound on the total interference, i.e., $\interfTotVar \le I_{\mathrm{split}}$. At first, this might seem counterintuitive, as we consider two possibly computationally challenging problems instead of one. However, as we show in Section~\ref{sec:Results}, by carefully choosing $\interference_{\mathrm{split}}$, this approach is very effective in practice. Furthermore, the following easy observation holds:

\begin{lemma}
\label{lem:split}
    An optimal solution $S_{\ge}^*$ to the problem $P_{\ge}$ is also optimal to QSTPI if $P_{\le}$ has no solution $\overline{S}_{\le}$ with $c(\overline{S}_{\le}) < c(S_{\ge}^*)$.
\end{lemma}

\begin{proof}
    Given an optimal solution $S_{\ge}^*$ to the problem $P_{\ge}$ and let $S^*$ be the optimal solution to QSTPI. If $c(S^*) < c(S_{\ge}^*)$, then $S^*$ is not feasible for $P_{\ge}$, but must be feasible for $P_{\le}$.
\end{proof}

\subsection{A valid lower bound for the total interference}\label{sec:MinI_lowerbound}
There are multiple ways to choose the split value $\interference_{\mathrm{split}}$, e.g., as a fixed parameter or based on some heuristic. Let us consider the following minimization problem:
\begin{align}
    \mathrm{MinI}:\quad\min\quad & \interfTotVar
    & \label{eq:MinInterf}\\
    \text{s.t.}\quad & \sum_{i \in \potTerminals} q_{i} y_i - \interfTotVar\ge Q & \label{eq:MinInterfQuota}\\
    & \interfTotVar \ge \sum_{i \in \potTerminals} \sum_{j>i \in \potTerminals} (I_{ij} + I_{ji}) y_j y_i
    \label{eq:MinInterfInterf}\\
    & y_i + y_j \le 1 & \forall i,j \in \potTerminals, i \neq j, d_{ij} < D_{\mathrm{min}}\label{eq:MinInterfMinDist}\\
    &y_i \in \lbrace0,1\rbrace & \forall i \in \potTerminals \\
    &\interfTotVar \in \mathbb{R}_{\ge 0},&\label{eq:MinInterfEnd}
\end{align}
i.e., the problem is to choose a subset of potential terminals that minimizes the total interference \eqref{eq:MinInterf}, while fulfilling the desired quota \eqref{eq:MinInterfQuota}, considering the total interference \eqref{eq:MinInterfInterf}, and respecting the minimal distance constraints \eqref{eq:MinInterfMinDist}. This problem can be seen as a version of the WFLO. As the total interference is minimized in $\mathrm{MinI}$, there will be no solution for QSTPI with a lower interference than in an optimal solution of $\mathrm{MinI}$. 

As splitting into up- and downproblem is pointless when $\interference_{\mathrm{split}} = 0$, we will now give a necessary and sufficient criterion when MinI has an optimal solution of value zero.
Let us consider the undirected graph $\graph_p = (\potTerminals, E_p)$, where for each pair $(i, j)$ of distinct vertices, there is an edge $(i, j)\in E_p$ if and only if $\interference_{ij} + \interference_{ji} > 0$ or $d_{ij} < D_{\mathrm{min}}$. A subset of vertices $S\subseteq\potTerminals$ is called an \emph{independent set} in $G_p$ iff there is no edge in $E_p$ that connects any two vertices in $S$. A \emph{maximum independent set} (MIS) of graph $\graph_p$ is an independent set that maximizes $\vertexProfits(S)$.
\begin{theorem}\label{thm:mis}
    Let $\quota_{\mathrm{MIS}}^*$ be the quota of an MIS in $\graph_p$ and $\interfTotVar^*$ an optimal solution to $\mathrm{MinI}$. Then
    \begin{align}
        \quota \le\quota_{\mathrm{MIS}}^* \Leftrightarrow \interfTotVar^* = 0
    \end{align}
\end{theorem}

\begin{proof}
    The total interference of $\mathrm{MinI}$ is bounded by zero from below. 
    Let $S$ be an MIS and set $y_i \coloneqq 1$ if and only if $i \in S$. By definition, $S$ does not contain any pair $(i, j)$ of vertices with $I_{ij} + I_{ji} > 0$ or $d_{ij} < D_{\mathrm{min}}$, so that \eqref{eq:MinInterfInterf} and \eqref{eq:MinInterfMinDist} trivially hold. Moreover, \eqref{eq:MinInterfQuota} reads as $\quota_{\mathrm{MIS}}^* - \interfTotVar\geq \quota$. In particular, $\interfTotVar= 0$ is feasible and hence optimal for $\mathrm{MinI}$ if and only if $\quota_{\mathrm{MIS}}^*\ge \quota$.
\end{proof}

Rephrasing Theorem~\ref{thm:mis}, a necessary condition for a non-trivial total interference ($\interfTotVar > 0$) is that $\quota > \quota_{\mathrm{MIS}}^*$. Combining Lemma~\ref{lem:split} and Theorem~\ref{thm:mis}, we obtain:

\begin{corollary}\label{Cor:MINI}
    Let $I_{\mathrm{LB}}$ be a valid lower bound on the problem $\mathrm{MinI}$. An optimal solution $S_{\ge}^*$ to the problem $P_{\ge}$ is also optimal to $P$ if $I_{\mathrm{LB}}$ is chosen as split value $I_{\mathrm{split}}$.
\end{corollary}

That means, by finding a valid lower bound to $\mathrm{MinI}$ greater than zero, only the smaller upproblem $P_\ge$ has to be solved to solve the full QSTPI.

\section{Computational Study}\label{sec:CompStudy}
In the following, we describe the implementation and the setup of our computational study in Section~\ref{sec:implement} and Section~\ref{sec:setting}, respetively. Furthermore, we present the data used in this study in Section~\ref{sec:data}. We conclude this section by demonstrating the strength of our approach compared to the state-of-the-art general MIP solver \gurobi in Section~\ref{sec:Results}.

\subsection{Implementation} \label{sec:implement}

Based on the promising results shown in \textcite{pedersen2024}, we have integrated QSTPI-TRANS as in Section~\ref{sec:qstpTrans} into the Steiner tree software package \scipjack. The exponentially-many constraints \eqref{eq:Steiner cut trans} are separated using the sophisticated maximum-flow algorithm implemented by \textcite{rehfeldt2021a} which is based on the classical max-flow/min-cut theorem. Additionally to the separation procedure, the solution process is supported by a \textit{shortest-path heuristic} (SPH). 

In general, heuristics are utilized to find good and feasible solutions in short time to provide upper bounds on the exact solution of the problem. Within the context of STP-related problems, the SPH is the best-known heuristic and was introduced by \textcite{takahashi1980}. For this study, we adapt \scipjack's SPH implementation by \textcite{rehfeldt2021a} and already customized for the case of the QSTP by \textcite{pedersen2024}: Starting at the root $r$, the closest fixed terminal or potential terminal is inserted to the \emph{root component}. The interference between the newly inserted node and all nodes in the root component is calculated and added to the total interference. Then, the distance of all non-connected vertices to the root component is updated. This process is repeated until all fixed terminals are connected and the connected potential terminals minus the total interference fulfill the desired quota. Preliminary experiments have shown that, although the SPH finds good upper bounds at the root node, it rarely improves the primal bound during the branch-and-bound procedure. Therefore, we decided to run the SPH only at the root node, modifying the arc costs to $\overline{\cost}_\arc \coloneqq (1.0 - \hat{\arcVariable}_\arc) \cost_\arc$ for all $\arc \in \setArcs$ depending on the current LP-solution $\hat{x}$. Our solver will be benchmarked against a black-box MIP solver using the flow-based formulation presented in Section~\ref{sec:flow}.

\subsection{Computational setting}\label{sec:setting}

The QSTPI-TRANS is integrated into \scipjack using \scipVersion{8}{0}{1}\parencite{bestuzheva2023} with \cplexVersion{12}{10} \parencite{ibmilog2022} as LP solver. The flow-based MIP formulation is implemented using the \gurobi~\python-interface in \python~\oldstylenums{3.11.2}, and is solved with \gurobiVersion{11}{01} \parencite{gurobi1101}. The following settings are run:
\begin{itemize}
    \item \setting{FLOW}: the basic flow-based model, presented in \eqref{eq:obj_flow} -- \eqref{eq:var_flow};
    \item \setting{FLOW-Ilb}: as before, but introducing the lower bound $I_{\mathrm{LB}}$ described in Section~\ref{sec:Split} on the total induced interference;
    \item \setting{QSTPI}: the QSTPI-TRANS formulation presented in \eqref{eq:obj_trans_qstp} -- \eqref{eq:quota_cons_trans} and \eqref{eq:minD_trans} -- \eqref{eq:bilinear_var};
    \item \setting{QSTPI-Ilb}: as before, but introducing the lower bound $I_{\mathrm{LB}}$  described in Section~\ref{sec:Split} on the total induced interference;
    \item \setting{QSTPI-S-H$\alpha$}: as before, but a $\interference_{\mathrm{split}}$ is determined using the interference value from the initial shortest-path heuristic solution multiplied by $\alpha \in \lbrace 0.1, 0.5, 0.8\rbrace$;
    \item \setting{QSTPI-S-minI-$\tau$}: as before, but $\interference_{\mathrm{split}} \coloneqq \max(I_{\mathrm{LB}}, \interfTotVar^*)$, where $\interfTotVar^*$ corresponds to the best primal solution of the minimization problem MinI presented in \eqref{eq:MinInterf} -- \eqref{eq:MinInterfEnd} within a time limit of $\tau~\in~\lbrace300, 600, 1800\rbrace$ seconds using \scipVersion{8}{0}{1}.
\end{itemize}

All computational experiments are executed single-threaded in a non-exclusive mode on a cluster with \textit{Intel(R) Xeon(R) Gold 6342} CPUs running at 2.80 GHz, where four CPUs and 64 GB of RAM are reserved for each run. For a better comparison with our implementation in \scipjack, which is only single-threaded, we restrict the black-box solver \gurobi{} to one thread when discussing the results. However, further experiments have shown that under the given setup no performance improvement is gained when running \gurobi{} in 8-threaded mode. We set a time limit of six hours (21,600 s). In the cases of \setting{QSTPI-S-H$\alpha$} and \setting{QSTPI-S-minI-$\tau$}, the downproblem is only solved if the upproblem is solved to optimality. The downproblem is then run using the optimal objective value of the upproblem as an objective bound (see Lemma~\ref{lem:split}). 

\subsection{Data creation}\label{sec:data}

\begin{figure}
\centering
\begin{subfigure}[b]{.3\linewidth}
\includegraphics[width=\linewidth]{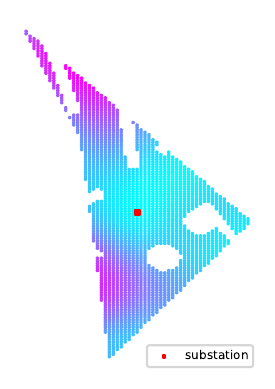}
\caption{Instance A}\label{fig:instance_A}
\end{subfigure}
\begin{subfigure}[b]{.3\linewidth}
\includegraphics[width=\linewidth]
{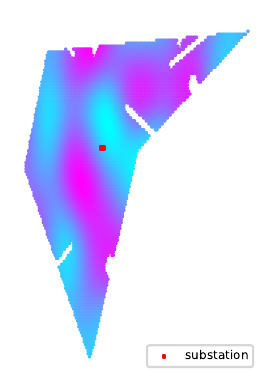}
\caption{Instance B}\label{fig:instance_B}
\end{subfigure}
\begin{subfigure}[b]{.3\linewidth}
\includegraphics[width=\linewidth]
{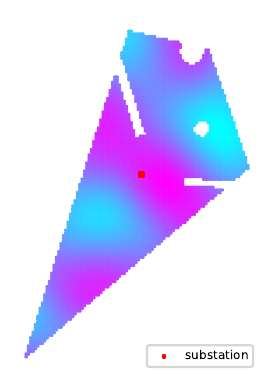}
\caption{Instance C}\label{fig:instance_C}
\end{subfigure}
\caption{Chosen representative wind farms by \cite{Cazzaro.2022}. Color of regions represents the costs of building a turbine from {\color{cyan}low} to {\color{magenta}high} costs.}
\label{fig:instances}
\end{figure}

As a benchmark set, we use three representative wind farm instances out of the dataset introduced by \textcite{Cazzaro.2022}, namely instances A, B, and C. Instance A and C are divided into two zones in the original data set of which we only consider the larger zone for this paper; see Fig.~\ref{fig:instances}. For each wind farm, we create 20 small-sized, medium-sized, and large-sized instances each by randomly sampling 100, 200, and 500 potential positions, respectively.

For instance A we set the quota equivalent to the annual power production without interference of 10 and 20 wind turbines with 15\,MW each based on the real wind measurements at Ten Noorden van de Waddeneiladen\footnote{reported by \url{https://english.rvo.nl/topics/offshore-wind-energy/noorden-waddeneilanden} Accessed on: 11/13/2025; Given in dataset from \textcite{Cazzaro.2022}} for all sampled instances. For instance B the quota is set to 10, 20, and 40 wind turbines for the small-sized instances and 10, 20, 40, and 60 for the medium- and large-sized instances in the same way as before. For instance C the quota of the small- and medium-sized instances is set to 10, 20 and 40, and for the large-sized ones to 10, 20, 40, and 60 wind turbines in the same way as before. Due to unlucky sampling, too many turbines might be too close to each other to satisfy the quota constraint. Removing these trivially infeasible instances, our complete benchmark set has 143 (instead of 160 instances) small-, 180 medium-, and 200 large-sized instances. As the focus lies on the number of potential turbine positions, no additional Steiner nodes are considered.

\begin{figure}
\centering
\begin{subfigure}[b]{.45\linewidth}
\centering
\includegraphics[width=.9\linewidth, trim=0 40 0 0, clip]{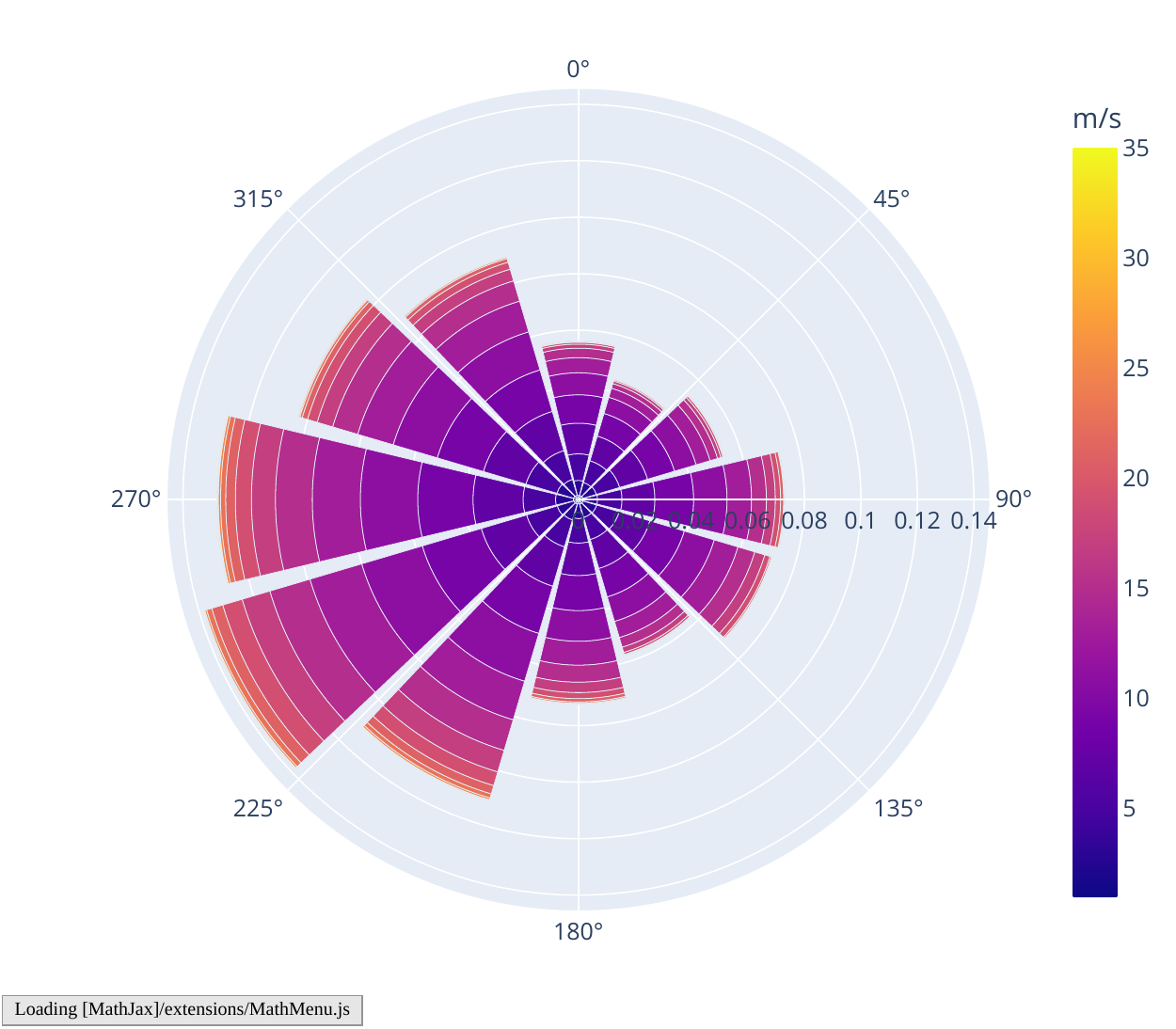}
\vspace{25pt}
\caption{Wind rose}\label{fig:windrose}
\end{subfigure}
\begin{subfigure}[b]{.45\linewidth}
\centering
\includegraphics[width=.7\linewidth]{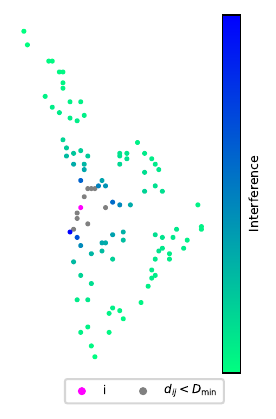}
\caption{Interference caused by turbine $i$}\label{fig:interf_scenario}
\end{subfigure}
\caption{Wind rose of Ten Noorden van de Waddeneilanden site (Panel~\ref{fig:windrose}); Interference caused by turbine $i$ in chosen wind scenario for a sampled instance of wind farm area A (Panel~\ref{fig:interf_scenario}); Interference calculated by the Jensen model \parencite{jensen1983}. Wind farm and wind data from \cite{Cazzaro.2022}.}
\label{fig:dataset}
\end{figure}

The interference is calculated by the sum of pairwise interferences using Jensen's model \parencite{jensen1983}, an analytical model for wind parks, as done in \textcite{fischetti2016proximity}, based on the same real wind measurements given in the dataset by \textcite{Cazzaro.2022}; see Fig.~\ref{fig:windrose}. 
While more accurate interference models exist in comparison to Jensen's model (e.g., \textcite{Gebraad2014Floris, Gebraad2014FlorisDyn, king2021, bay2023}), it still provides a useful approximation for early-stage wind farm layout planning \parencite{fischetti2016proximity, archer2018wakereview}. As our focus is on network design methodology, and we assume that the interference is given as an input to our problem, we stick with this computationally simpler model.
Fig.~\ref{fig:interf_scenario} shows the interference caused by a turbine under the considered wind scenario in an exemplary sample of instance A. No information regarding the position of the substation is given in the original dataset. Therefore, the substation is placed in the middle of each area, i.e., the coordinates of the substation are given by averaging the coordinates of all available positions. The dataset by \textcite{Cazzaro.2022} is used, which contains the cost to build a turbine at a specific position, which includes the foundation costs, as well as data for 15 MW turbines with a rotor diameter of 240 meters. The minimum distance required between two turbines is usually three to five times the rotor diameter; Following \textcite{fischetti2022}, we set the minimum distance $D_{\mathrm{min}}$ to five times the rotor diameter, i.e., 1200 meters. The cable costs are chosen in accordance to \textcite{Cazzaro.2020}. The authors provide cable costs for three types of cables: 430, 480, and 610 k€/km. Since we do not consider different cable types, a mean value of 504 k€/km is chosen.

\subsection{Results}\label{sec:Results}

\begin{figure}
    \centering
    \includegraphics[width=.8\linewidth]{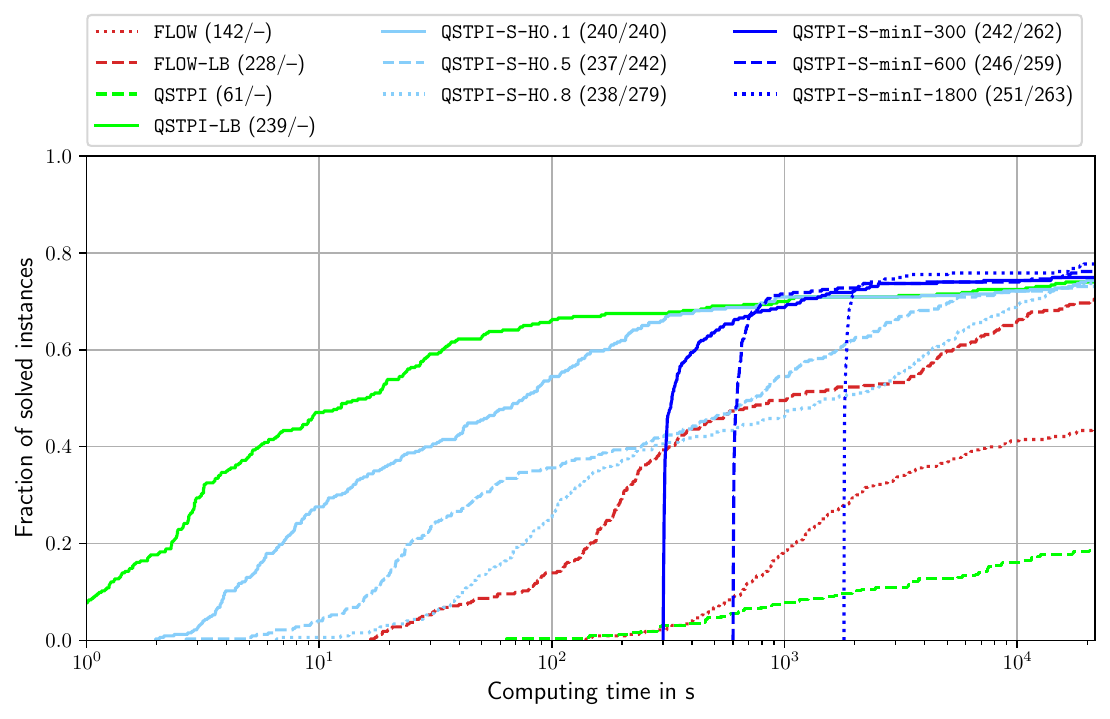}
    \caption{Cumulative percentage of instances solved to optimality over time for all settings for small- and medium-sized instances. The numbers in brackets state the number of instances which are solved to optimality completely and instances for which the upproblem is solved to optimality.}
    \label{fig:perfProf.offshore1}
\end{figure}

Fig.~\ref{fig:perfProf.offshore1} shows the cumulative percentage of instances solved to optimality over time for the complete problem for all settings. In the legend of the figure, the first number in the brackets reports the number of instances completely solved to optimality. The second number gives the number of instances, for which the upproblem was solved to optimality. In all solved instances, the global optimal solution lies within the upproblem. Furthermore, the total interference induces by the optimally chosen turbines is far away of the introduced splitting value $\interference_{\mathrm{split}}$, when considering all settings but \setting{QSTPI-S-H0.8}.

Whereas the simple QSTP without the interference and minimum distance constraints has outperformed the flow-based formulation by up-to two orders of magnitude before (see \textcite{pedersen2024}), this advantage vanishes complete considering \setting{QSTPI} and \setting{FLOW}. \setting{QSTPI} gives very poor performance only solving 61 instances within the time limit of six hours, which is not even half the number of instances solved by \setting{FLOW}. The performance is drastically increased by introducing the lower bound $I_{\mathrm{LB}}$ described in Section~\ref{sec:Split}. \setting{FLOW-Ilb} is not only more than a magnitude faster than \setting{QSTPI} and \setting{FLOW}, but also solves 167 and 86 more instances, respectively. Applying the same lower bound value on our approach in \setting{QSTPI-Ilb}, the performance is increased by another magnitude solving 239 of 323 instances to optimality. Remarkably, all these instances are solved within roughly one hour: In relation to that, \setting{FLOW-Ilb} only solves around 50\% of the instances in that time.

Considering \setting{QSTPI-S-H$\alpha$}, the best performance is achieved with $\alpha=0.1$, being around a magnitude faster than $\alpha=0.5$ and $\alpha=0.8$. \setting{QSTPI-S-H0.1} solves 240, \setting{QSTPI-S-H0.5} solves 237, and \setting{QSTPI-S-H0.8} solves 238 instances out of 323 instances to optimality. By using a higher split value the better is the performance in solving the number of instance for the upproblem, i.e., 240, 242, and 279 for $\alpha$ equal to 0.1, 0.5, and 0.8, respectively (see the second number in the brackets in the legend of Fig.~\ref{fig:perfProf.offshore1}). However, the downproblem becomes much harder to solve, and 0, 5, and 41 instances run into the time limit in the downproblem for $\alpha$ equal 0.1, 0.5, and 0.8, respectively.

For the settings \setting{QSTPI-S-minI-$\tau$} none of the starting problem MinI that determines the split value $I_{\mathrm{split}}$ is solved to optimality within the time limit of $\tau$. Therefore, the minimum time to solve an instances using \setting{QSTPI-S-minI-$\tau$} is greater than $\tau$. However, the overall performance is similar to \setting{QSTPI-Ilb} with the advantage of solving three, seven, and 12 more instances for $\tau$ equals 300, 600, and 1800 second, respectively. 

\begin{figure}
\centering
\begin{subfigure}[b]{.45\linewidth}
\includegraphics[width=\linewidth]{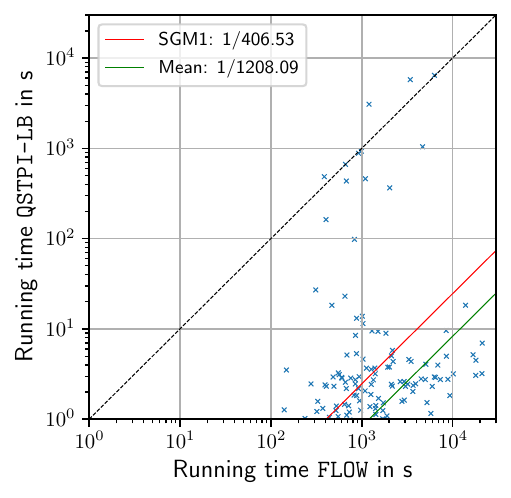}
\caption{\setting{QSTPI-Ilb} vs. \setting{FLOW}: 141 instances solved by both settings.}\label{fig:speedup-qstpIlb-grbbase}
\end{subfigure}
\begin{subfigure}[b]{.45\linewidth}
\includegraphics[width=\linewidth]{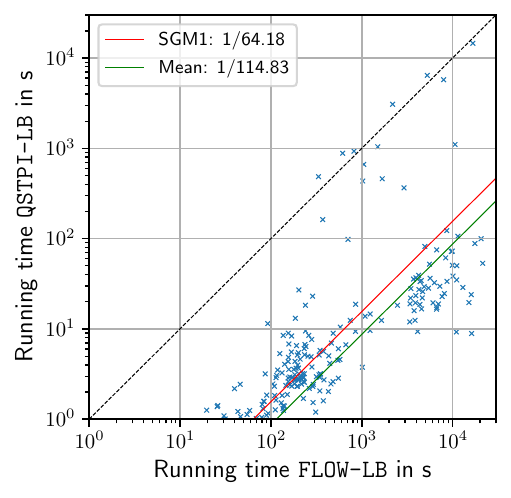}
\caption{\setting{QSTPI-Ilb} vs. \setting{FLOW-Ilb}: 227 instances solved by both settings.}\label{fig:speedup-qstpIlb-grbIlb}
\end{subfigure}
\caption{Comparison of running times between \setting{QSTPI-Ilb} and \setting{FLOW} (Panel~\ref{fig:speedup-qstpIlb-grbbase}) and of running times between \setting{QSTPI-Ilb} and \setting{FLOW-Ilb} (Panel~\ref{fig:speedup-qstpIlb-grbIlb}) of all instances solved by both settings, respectively. The dotted line shows the break even line; The \textit{SGM1} line represents the shifted geometric mean of the speedup with a shift of one second; the \textit{Mean} line represents the arithmetic mean of the speedup.}
\label{fig:speedup-qstp-grb}
\end{figure}

The performance profiles in Fig.~\ref{fig:perfProf.offshore1} have already shown the advantages of our novel approach for the IWFLCR problem. In the following, the setting \setting{QSTPI-Ilb} is first compared instance-wise to \setting{FLOW} and then to \setting{FLOW-Ilb} in Fig.~\ref{fig:speedup-qstpIlb-grbbase} and Fig.~\ref{fig:speedup-qstpIlb-grbIlb}. For the pairwise comparison we only consider instances solved to optimality by both settings. In the legend, additional to the arithmetic mean of the speedup, the shifted geometric mean with a shift of 1 (SGM1) is reported. The \textit{shifted geometric mean} \parencite{achterberg2007Phd} has become a standard measure in discrete optimization \parencite{mittelmann2020} and is less biased towards harder instances. By using \setting{QSTPI-Ilb} instead of \setting{FLOW} an average speedup of 1208.9 and a SGM1 of 406.53 is achieved, where \setting{QSTPI-Ilb} is faster on all but 5 instances, see Fig.~\ref{fig:speedup-qstpIlb-grbbase}. \setting{QSTPI-Ilb} is 114.83 times faster on average and 64.18 times faster considering the SGM1 compared to \setting{FLOW-Ilb}; see~\ref{fig:speedup-qstpIlb-grbIlb}. Again, \setting{QSTPI-Ilb} is faster on all but 5 instances.

\begin{figure}
    \centering
    \includegraphics[width=.8\linewidth]{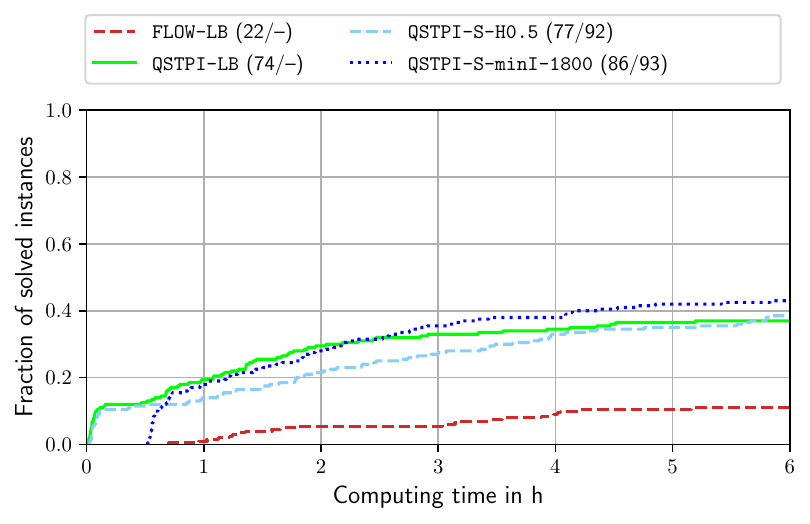}
    \caption{Cumulative percentage of instances solved to optimality over time for promising settings for large-sized instances. The numbers in brackets state the number of instances which are solved to optimality completely and instances for which the upproblem is solved to optimality.}
    \label{fig:perfProf.offshore500}
\end{figure}

The benefit of our novel approach becomes even more prominent when we are looking at the performance for the large-sized instances in Fig.~\ref{fig:perfProf.offshore500}. For these instances, only settings \setting{FLOW-Ilb}, \setting{QSTPI-Ilb}, \setting{QSTPI-S-H0.5}, and  \setting{QSTPI-S-minI-1800} are run, as they are expected to perform best for the large-sized instances. The setting \setting{FLOW-Ilb} only solves 22 instances out of 200 to optimality. The performance of the QSTPI related settings is relatively similar. However, the best setting is \setting{QSTPI-S-minI-1800} solving 86 out of 200 instances to optimality and seven more instances in the upproblem. The setting \setting{QSTPI-S-H0.5} solves 77 instances and with that three more instances than \setting{QSTPI-Ilb}.

\section{A capacitated approach for the QSTPI}
\label{sec:hop-qstpi}

\begin{figure}
    \centering
    \includegraphics[width=0.5\linewidth]{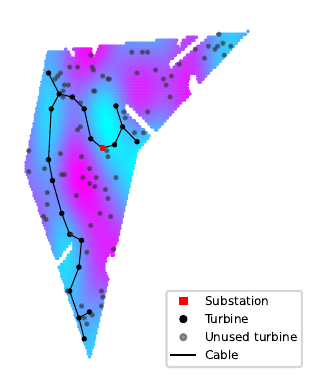}
    \caption{Exemplary unbalanced solution for the QSTPI connecting 16 turbines to the substation along a single. If, e.g., as in \cite{cazzaro2023}, a cable can only connect at most six 15\,MW turbines, this is not a practically feasible solution.}
    \label{fig:exemplary_qstpi_base}
\end{figure}

In the previous section, the advantages of formulating the IWLCR problem as Steiner-tree-based mixed integer program instead of a flow-based approach have become clear: We solve larger instances in shorter computational time. On the downside, our approach might produce solutions, in which a higher number of turbines are connected via a single cable with the substation, violating practically required cable capacity limits; see Fig.~\ref{fig:exemplary_qstpi_base}. These restrictions can be dealt with in an a-posteriori step, and we will do so in Section~\ref{sec:Res_Sequentiel_vs_combined}. However, we can also integrate cable capacity constraints into our previous model.

\subsection{Hop-constrained QSTPI}

Various topologies are possible for the cable routing in a wind farm, see, e.g., \textcite{Fischetti.2018} and \textcite{perezrua2019}. A common requirement in industrial application is a radial string topology, in which turbines are connected in (balanced) strings outgoing of the substation, see, e.g., \textcite{cazzaro2022, cazzaro2023}. Exploiting this fact, we model cable capacities by limiting the length of each string instead of imposing capacity bounds on each arc as usually done in flow-based network planning approaches. This limit is referred to as \textit{hop-limit}, i.e., the limit of arcs one is allowed to traverse (\textit{hop}) to reach a chosen turbine. For more details on general hop-constrained Steiner tree problems, the reader is referred to \textcite{costa2009hop, gouveia2011, sinnl2016node}.
We call our extension by hop limits the \emph{hop-constrained QSTPI}.

\subsubsection{Layered graph approach}

To model the hop-constrained QSTPI, we follow the layered graph approach proposed by \textcite{sinnl2016node}. We consider the transformed directed graph $G = (V, A)$ as described in Section~\ref{sec:qstpTrans}. We distinguish between three pairwise disjoint types of $V$: The potential terminals $i \in T_p$ (turbines) with corresponding fixed terminals $i' \in T_f$, and the artificial root (substation) $r$. We assume that no further fixed terminals exist, as our aim is to build a radial layout, and we also disregard additional Steiner nodes. Given a \textit{hop-limit} $H$, we create a layered graph $G_L = (V_L,A_L)$ with $H+1$ layers.  The set of nodes is defined as $V_L \coloneqq V_0 \cup \ldots\cup V_H$, where
\begin{equation}
\begin{aligned}
    V_0 &\coloneqq\lbrace r\rbrace \cup \fixTerminals,\\
    V_h &\coloneqq \{ i_h \mid i \in T_p \}, \quad h \in \lbrace 1, \dots, H\rbrace,
\end{aligned}
\end{equation}
i.e., layer zero only contains the root node and the fixed terminals, while all other layers contain a copy of each original potential terminal. The node $i_h$ refers to the node $i$ in layer $h$.
The set of arcs is defined as $A_L\coloneqq A_0\cup\ldots \cup A_{H}$, where
\begin{equation}
    \begin{aligned}
        A_0 &\coloneqq \lbrace (r,i') \mid (r,i') \in A, i' \in T_f\rbrace ,\\
        A_1 &\coloneqq \lbrace (i_h, i^\prime) \mid (i, i^\prime) \in A, i\in T_p, i^\prime \in T_f, h \in \lbrace1,\dots, H\rbrace\rbrace,\\
        A_h &\coloneqq \lbrace(i_{h-1}, j_{h}) \mid (i,j) \in A, i, j \in T_p\rbrace, & h \in \{2, \dots, H\}.
    \end{aligned}
\end{equation}
Figure~\ref{fig:GraphTransHop} shows the transformation from the original graph to the layered graph with a hop-limit of $H=3$. 

\begin{figure}
\centering
    \begin{subfigure}{0.4\linewidth}
        \centering
        \begin{tikzpicture}[scale=1.5, every node/.style={scale=1.}, rotate=90, yscale=-1]
            \tikzstyle{terminal} = [draw, fill, thick, minimum size=0.4, inner sep=2.6pt, color=red];
            \tikzstyle{pot_terminal} = [draw=blue, pattern=north west lines, pattern color=blue, thick, minimum size=0.4, inner sep=2.6pt];
            
            \def\y{16}
            \node[terminal] (r) at (.75+\y, 3.5) {} node[above of=r,yshift=-20pt]{$r$};
            \node[pot_terminal] (pt1) at(0. + \y, 2.5) {} 
            node[above of=pt1, xshift=-5pt, yshift=-20pt]{\scriptsize $1$};
            \node[pot_terminal] (pt2) at(0. + \y, 4.5) {} 
            node[above of=pt2, xshift=-5pt, yshift=-20pt]{\scriptsize $2$};

            \node[pot_terminal] (pt3) at(-1. + \y, 2.5) {} 
            node[above of=pt3, xshift=-5pt, yshift=-20pt]{\scriptsize $3$};

            \node[pot_terminal] (pt4) at(-1. + \y, 4.5) {} 
            node[above of=pt4, xshift=-5pt, yshift=-20pt]{\scriptsize $4$};

            \draw[-,thick] (r) -- (pt1);
            \draw[-,thick] (r) -- (pt2);

            \draw[-,thick] (pt1) -- (pt3);
            \draw[-,thick] (pt1) -- (pt4);

            \draw[-,thick] (pt2) -- (pt3);
            \draw[-,thick] (pt2) -- (pt4);

            \draw[-,thick] (pt3) -- (pt4);
        \end{tikzpicture}
        \caption{Original instance.}
    \end{subfigure}
    \begin{subfigure}{0.58\linewidth}
    \centering
       \begin{tikzpicture}[scale=1.5, every node/.style={scale=1.}, rotate=90, yscale=-1]
            \tikzstyle{terminal} = [draw, fill, thick, minimum size=0.4, inner sep=2.6pt, color=red];
            \tikzstyle{pot_terminal} = [draw=blue, pattern=north west lines, pattern color=blue, thick, minimum size=0.4, inner sep=2.6pt];
            
            \def\y{16}
            \node[terminal] (r) at (.75+\y, 3.5) {} node[above of=r,yshift=-20pt]{$r$};
            
            \node[terminal] (t1) at(0.75 + \y, 2.5) {} 
            node[above of=t1, xshift=-5pt, yshift=-20pt]{\scriptsize $1^\prime$};
            \node[terminal] (t2) at(0.75 + \y, 4.5) {} 
            node[above of=t2, xshift=-5pt, yshift=-20pt]{\scriptsize $2^\prime$};
            \node[terminal] (t3) at(0.75 + \y, 1.5) {} 
            node[above of=t3, xshift=-5pt, yshift=-20pt]{\scriptsize $3^\prime$};
            \node[terminal] (t4) at(0.75 + \y, 5.5) {} 
            node[above of=t4, xshift=-5pt, yshift=-20pt]{\scriptsize $4^\prime$};
            
            \node[pot_terminal] (pt1) at(0. + \y, 2.5) {} 
            node[above of=pt1, xshift=-5pt, yshift=-20pt]{\scriptsize $1$};
            \node[pot_terminal] (pt2) at(0. + \y, 4.5) {} 
            node[above of=pt2, xshift=-5pt, yshift=-20pt]{\scriptsize $2$};
            \node[pot_terminal] (pt3) at(-1. + \y, 2.5) {} 
            node[above of=pt3, xshift=-5pt, yshift=-20pt]{\scriptsize $3$};
            \node[pot_terminal] (pt4) at(-1. + \y, 4.5) {} 
            node[above of=pt4, xshift=-5pt, yshift=-20pt]{\scriptsize $4$};

            \draw[-,thick] (r) -- (pt1);
            \draw[-,thick] (r) -- (pt2);

            \draw[-,thick] (pt1) -- (pt3);
            \draw[-,thick] (pt1) -- (pt4);

            \draw[-,thick] (pt2) -- (pt3);
            \draw[-,thick] (pt2) -- (pt4);

            \draw[-,thick] (pt3) -- (pt4);

            \draw[thick, -latex](r) edge[bend left=0]  (t1);
            \draw[thick, -latex](r) edge[bend left=0]  (t2);
            \draw[thick, -latex](r) edge[bend left=30]  (t3);
            \draw[thick, -latex](r) edge[bend left=-30]  (t4);

            \draw[thick, -latex](pt1) edge[bend left=0]  (t1);
            \draw[thick, -latex](pt2) edge[bend left=0]  (t2);
            \draw[thick, -latex](pt3) edge[bend left=-20]  (t3);
            \draw[thick, -latex](pt4) edge[bend left=20]  (t4);
        \end{tikzpicture}
        \caption{Transformed instance.}
    \end{subfigure}
    \begin{subfigure}{\linewidth}
    \centering
       \begin{tikzpicture}[scale=1.5, every node/.style={scale=1.}, rotate=90, yscale=-1]
            \tikzstyle{terminal} = [draw, fill, thick, minimum size=0.4, inner sep=2.6pt, color=red];
            
            \tikzstyle{pot_terminal} = [draw=blue, pattern=north west lines, pattern color=blue, thick, minimum size=0.4, inner sep=2.6pt];
            
            \def\y{16}
            \node[terminal] (r) at (.75+\y, 3.5) {} node[above of=r,yshift=-20pt]{$r$};
            
            \node[terminal] (t1) at(0.75 + \y, 2.25) {} 
            node[above of=t1, xshift=0pt, yshift=-22pt]{\scriptsize ${1^\prime}$};
            \node[terminal] (t2) at(0.75 + \y, 4.75) {} 
            node[above of=t2, xshift=0pt, yshift=-22pt]{\scriptsize ${2^\prime}$};
            \node[terminal] (t3) at(0.75 + \y, 1.5) {} 
            node[above of=t3, xshift=0pt, yshift=-22pt]{\scriptsize ${3^\prime}$};
            \node[terminal] (t4) at(0.75 + \y, 5.5) {} 
            node[above of=t4, xshift=0pt, yshift=-22pt]{\scriptsize ${4^\prime}$};

            \node[ellipse, draw=cyan , fill=cyan, fill opacity=0.05,
                fit=(t1)(t2)(t3)(t4), inner ysep=8pt]  (group) {};
            \node[anchor=west] at (group.west) {$L_0$};
            
            \node[pot_terminal] (pt11) at(-0.12 + \y, 2.5) {} 
            node[above of=pt11, xshift=-5pt, yshift=-20pt]{\scriptsize $1_1$};
            \node[pot_terminal] (pt21) at(-0.12 + \y, 4.5) {} 
            node[above of=pt21, xshift=-5pt, yshift=-20pt]{\scriptsize $2_1$};

            \node[ellipse, draw=magenta, fill=magenta, fill opacity=0.05,
                fit=(pt11)(pt21), inner ysep=8pt, inner xsep=47pt] (group) {};
            \node[anchor=west] at (group.west) {$L_1$};
                
            \node[pot_terminal] (pt32) at(-1. + \y, 2.5) {} 
            node[above of=pt32, xshift=-5pt, yshift=-20pt]{\scriptsize $3_2$};
            \node[pot_terminal] (pt42) at(-1. + \y, 4.5) {} 
            node[above of=pt42, xshift=-5pt, yshift=-20pt]{\scriptsize $4_2$};
            \node[ellipse, draw=cyan , fill=cyan, fill opacity=0.05,
                fit=(pt32)(pt42), inner ysep=8pt, inner xsep=47pt] (group) {};
            \node[anchor=west] at (group.west) {$L_2$};

            \node[pot_terminal] (pt33) at(-1.9 + \y, 2.5) {} 
            node[below of=pt33, xshift=7pt, yshift=23pt]{\scriptsize $3_3$};
            \node[pot_terminal] (pt43) at(-1.9 + \y, 4.5) {} 
            node[below of=pt43, xshift=7pt, yshift=23pt]{\scriptsize $4_3$};
            \node[pot_terminal] (pt13) at(-1.9 + \y, 1.5) {} 
            node[below of=pt13, xshift=7pt, yshift=23pt]{\scriptsize $1_3$};
            \node[pot_terminal] (pt23) at(-1.9 + \y, 5.5) {} 
            node[below of=pt23, xshift=7pt, yshift=23pt]{\scriptsize $2_3$};
            
            \node[ellipse, draw=magenta , fill=magenta, fill opacity=0.05,
                fit=(pt13)(pt23)(pt33)(pt43), inner ysep=8pt] (group) {};
            \node[anchor=west] at (group.west) {$L_3$};

            \draw[-latex] (r) -- (pt11);
            \draw[-latex] (r) -- (pt21);

            \draw[-latex] (pt11) -- (pt32);
            \draw[-latex] (pt11) -- (pt42);

            \draw[-latex] (pt21) -- (pt32);
            \draw[-latex] (pt21) -- (pt42);

            \draw[-latex] (pt32) -- (pt23);
            \draw[-latex] (pt32) -- (pt43);
            \draw[-latex] (pt32) -- (pt13);

            \draw[-latex] (pt42) -- (pt33);
            \draw[-latex] (pt42) -- (pt13);
            \draw[-latex] (pt42) -- (pt23);

            \draw[ -latex](r) edge[bend left=0]  (t1);
            \draw[ -latex](r) edge[bend left=0]  (t2);
            \draw[ -latex](r) edge[bend left=30]  (t3);
            \draw[ -latex](r) edge[bend left=-30]  (t4);

            \draw[ -latex](pt11) edge[bend left=0]  (t1);
            \draw[ -latex](pt21) edge[bend left=0]  (t2);
            \draw[ -latex](pt32) edge[bend left=-20]  (t3);
            \draw[ -latex](pt42) edge[bend left=20]  (t4);

            \draw[-latex](pt13) edge[bend left=0]  (t1);
            \draw[-latex](pt23) edge[bend left=0]  (t2);

            \draw[-latex](pt33) edge[bend left=-20]  (t3);
            \draw[-latex](pt43) edge[bend left=20]  (t4);
        \end{tikzpicture}
        \caption{Layered instance with hop-limit $H=3$.}
    \end{subfigure}
    \caption{Transformation for the hop-constrained QSTPI. We add first the artificial root $r$, and for each potential terminal $i \in \potTerminals$ (blue squares) a fixed terminal (red square) $i^\prime \in T_f$ as before to obtain the transformed graph $G$. We then build the layered graph $G_L$ with $H = 3$ layers. We omit nodes with no ingoing arcs, as they cannot be reached by a directed path from $r$.} 
    \label{fig:GraphTransHop}
\end{figure}

The rationale behind the layered graph construction is the following observation,
that allows to translate between certain Steiner trees in the transformed graph $G$ and hop-constrained radial-layout Steiner trees in $G_L$.
\begin{lemma}
    There is a one-to-one correspondence between the following sets:
    \begin{itemize}
        \item the arborescences in $G_L$ that are rooted at $r$, contain $T_f$ and at most one node $i^h$ for each $i \in T_p$, and have the property that each $i^h \in V_1 \cup \dots \cup V_{H-1}$ has at most one outgoing arc to layer $h+1$.
        \item the arborescences in $G$ that are rooted at $r$, contain $T_f$, and have the property that for all $i \in T_p$, the number of outgoing arcs is at most one, and that the length of the unique $r$-$i$-path is at most $H$.
    \end{itemize}
\end{lemma}
\begin{proof}
    If $S_L$ is an arborescence in $G_L$ with the above properties, then $S_L$ is necessarily a union of directed $r$-$i^h$-paths for some $i^h \in V_h$ and directed $r$-$i'$-paths for some $i' \in T_p$. By construction of $A_L$, arcs to nodes in $V_1 \cup \dots \cup V_H$ always go to the next layer, so that all nodes in $V_1 \cup \dots \cup V_{H-1}$ must be reached from $r$ using at most $H$ arcs. As for all $i \in T_p$, $S_L$ contains at most one of the $i^h$, so that, replacing $i^h$ by $i$, we obtain an arborescence in $S$ with the desired properties. For the converse, we only need to track the layer count.
\end{proof}

\subsubsection{Model formulation}

Before describing the model formulation in detail, we define an additional set of variables $X$ for convenience. For each arc $(i,j) \in A$ with $i,j \in T_p$, let
\begin{align}
    X_{ij} \coloneqq \sum_{h=2}^H\arcVar_{i_{h-1}j_{h}},
\end{align}
and for each arc $(i,i^\prime)\in A$ with $i\in T_p$ and $i' \in T_f$, let
\begin{align}
    X_{ii^\prime} \coloneqq \sum_{h=1}^H\arcVar_{i_h i^\prime}.
\end{align}
These variables link an arc of the original graph with the arcs in the layered graph: If arc $(i,j)\in A$ is chosen for the original graph ($X_{ij} = 1$), it must connect two consecutive layers. If the fixed terminal $i^\prime \in T_f$ is not connected via the root node $r$ ($x_{ri^\prime}=0$ and $X_{ii^\prime} = 1$), it must be connected via any layer node $i_h$ with $h\in\lbrace1,\dots, H\rbrace$. The hop-constrained QSTPI is formulated as:

\begin{align}
    \text{QSTPI-HOP:}&&\nonumber\\
    \quad\min \quad&\cost^T \arcVar &\label{eq:obj_qhop}\\
    \text{s.t.} \quad
    & \arcVar(\delta^-(W)) \ge 1& \forall W \subset \setVertices_L, r\notin W, |W \cap \fixTerminals| \ge 1\label{eq:Steiner_qhop}\\
    & \sum_{i \in \potTerminals} \vertexProfitsNumbered{v} \arcVar_{ri^{\prime}} + \interfTotVar \le \sum_{i \in \potTerminals} \vertexProfitsNumbered{i} - \quota & \label{eq:quota_cons_qhop}\\
    &\interfTotVar \ge \sum_{i \in \potTerminals} \sum_{j \neq i \in \potTerminals} \interference_{ij} X_{ii^{\prime}} X_{jj^{\prime}}&\label{eq:interf_tot_qhop}\\
    & X_{ii^{\prime}} + X_{jj^{\prime}} \le 1 & \forall i,j \in \potTerminals, i \neq j, d_{ij} < D_{\mathrm{min}}\label{eq:minD_qhop}\\
    & \sum_{(i,j)\in A, j\in T_p} X_{ij} \le 1 & i\in T_p\label{eq:Deg2_qhop}\\
    & X_{ij} = \sum_{h=2}^H\arcVar_{i_{h-1}j_{h}}& \forall (i,j) \in A, i,j \in T_p\label{eq:Xvw_qhop}\\
    & X_{ii^\prime} = \sum_{h=1}^H\arcVar_{i_h i^\prime}& \forall i\in T_p\label{eq:Xvv_qhop}\\
    & \arcVar_{ij} \in \lbrace 0,1 \rbrace &\forall (i,j) \in \setArcs_L\label{eq:var_qhop}
\end{align}

Constraint \eqref{eq:Steiner_qhop} -- \eqref{eq:minD_qhop} are inherited from the QSTPI-TRANS (see Section~\ref{sec:qstpTrans}), i.e., each terminal $i^\prime\in T_f$ must be connected to the root node $r$ in \eqref{eq:Steiner_qhop}, the quota must be fulfilled in \eqref{eq:quota_cons_qhop}, the total interference is determined by the chosen turbines in \eqref{eq:interf_tot_qhop}, and only one of two turbines that are too close to each other can be chosen in \eqref{eq:minD_qhop}.
Constraint \eqref{eq:Deg2_qhop} ensures the radial topology by limiting the number of outgoing arcs of any $i^h$ to one. Note that \eqref{eq:Deg2_qhop} technically still allows that the last node in a path of the radial layout belongs to a potential terminal that has been visited before, thus not contributing to the quota, but possibly increasing the interference and the total cost, so that this phenomenon will not occur in an optimal solution. Similarly, $X_{ii'} > 1$ is not explicitly excluded in our model, but we can set $X_{ii'} = 1$ without further consequences in a post-processing step.
Constraints \eqref{eq:Xvw_qhop} and \eqref{eq:Xvv_qhop} link the arc variables of the layered graph to the arcs variables in the original graph. Note that \eqref{eq:Deg2_qhop} implies $X_{ij} \leq 1$.

\subsubsection{Additionally valid inequalities for string topologies}

\paragraph{Subtour elimination.}

The following generalized subtour elimination constraint of size two is valid \parencite{sinnl2016node}:
\begin{equation}
\begin{aligned}
    & X_{ij} + X_{ji} \le  X_{ii^\prime} & \forall (i, j) \in A, i\in T_p.\label{eq:SEC2_qhop}
\end{aligned}
\end{equation}
Constraint \eqref{eq:SEC2_qhop} states that whenver two nodes $i,j \in T_p$ are connected, then the fixed terminal $i^\prime$ must be connected from some layer.

Using the fact that there is an optimal solution, where $i^\prime\in T_f$ is either connected directly via the root  ($x_{ri^\prime} = 1$) or via a layer node ($X_{ii^\prime} = 1)$, we can use 
    $x_{ri^\prime} + X_{ii^\prime} = 1$ to rewrite constraint \eqref{eq:SEC2_qhop} as 
\begin{align}
    & X_{ij} + X_{ji} \le 1 - \arcVar_{ri^\prime} & \forall (i, j) \in A, i, j\in T_p\label{eq:rSEC2_qhop},
\end{align}
which empirically yields much better performance.

\paragraph{Linking the layers.}

So far, our integer program has considerable symmetry with respect to the layers. We will now exploit that we are aiming at a radial topology. For every layer node $i_h$, let 
\begin{align*}
    y_{i_h} \coloneqq \sum_{(j_{h-1}, i_h) \in A_h} x_{j_{h-1} i_h}
\end{align*}
describe whether we choose the potential terminal $i \in T_p$ in layer $h$. 
As we build an arborescence, each node cannot have more than one incoming arc, thus, $y_v^h$ describes if a node $v$ appears in layer $h$. 
Considering a node $v^h$ in layer $h\in\lbrace1,\dots, H-1\rbrace$, an arc only can leave the node if it appears in that layer. Thus, the following holds
\begin{align}
    &\sum_{(i_h, j_{h+1}) \in A_h} \arcVar_{i_h j_{h+1}} \le y_{i_h} & \forall i_h \in V_h, h \in \lbrace1,\dots,H-1\rbrace.\label{eq:qhop_inflow}
\end{align}
Additionally, in a radial layout, a layer $h\in\lbrace2,\dots,H\rbrace$ cannot contain more nodes than its previous layer $h-1$, as, otherwise, the solution would contain a node in layer $h-1$ having more than one arc to layer $h$:
\begin{align}
&\sum_{i_h\in V_h} y_{i_h} \le \sum_{i_{h-1}\in V_{h-1}} y_{i_{h-1}} &\forall h \in \lbrace2, \dots, H\rbrace.\label{eq:layerConsecutive}
\end{align}

\paragraph{Bounding the number of chosen turbines.}
Considering the quota constraint
\begin{align}
    \sum_{i\in T_p}q_i y_i - I_{\mathrm{tot}} \ge Q,
\end{align}
cf.\ \eqref{eq:quota_flow}, we can impose a valid upper bound on the number of chosen turbines. Let $q^{\mathrm{min}} \coloneqq\left( q_1, \dots, q_{|T_p|}\right)$ be sorted in non-decreasing order, $q_1\le q_2 \le \dots \le q_{|T_p|}$, and let $q^{\mathrm{min}}_i$ be the $i$-th entry of $q^{\mathrm{min}}$. Let $Q^{\mathrm{min}}(k)\coloneqq\sum_{i=1}^k q^{\mathrm{min}}_i$ be the minimum amount of quota collected by $k$ turbines. 
Given the interference matrix $I$, let $I^{\mathrm{max}}_i(k)$ denote the sum of the $k-1$ maximum $I_{ij}$ values of row $i$. We can then compute the worst case induced total interference $I^{\mathrm{max}}(k) \coloneqq k\max_i I^{\mathrm{max}}_i(k)$. A valid upper bound on the number of chosen turbines is given by the smallest $k^{\mathrm{ub}}$ for which 
\begin{align}
    Q^{\mathrm{min}}(k^{\mathrm{ub}}) - I^{\mathrm{max}}(k^{\mathrm{ub}})&\ge Q \\
    \text{and}\quad Q^{\mathrm{min}}(k^{\mathrm{ub}}-1) - I^{\mathrm{max}}(k^{\mathrm{ub}}-1) &< Q.
\end{align}
This upper bound in mind, given $H$ layers and the fact that no layer can contain more nodes than its predecessor, see \eqref{eq:layerConsecutive}, the following holds:
\begin{align}
    \sum_{(i_{H-1}, j_H)\in A_H} \arcVar_{i_{H-1} j_H} \le \left\lfloor \frac{k^{\mathrm{ub}}}{H}\right\rfloor,&\label{eq:arcLastlayerLim}\\
    \sum_{(i_H, i^\prime) \in A_0, i_H \in V_H} x_{i_{H} i^\prime} \le \left\lfloor \frac{k^{\mathrm{ub}}}{H}\right\rfloor.&\label{eq:nodeLastlayerLim}
\end{align}
Inequality \eqref{eq:arcLastlayerLim} limits the number of arcs connected to the last layer and \eqref{eq:nodeLastlayerLim} limits the nodes that connect to the terminals from the last layer. For instance, if $k^{\mathrm{ub}} = 22$ and $H = 6$, then at most $\left\lfloor \frac{22}{6}\right\rfloor = \left\lfloor 3.667\right\rfloor = 3$ turbines appear in the last layer, and, thus, only three arcs will enter the last layer.

\subsection{Computational Study}
We integrated the QSTPI-HOP model \eqref{eq:obj_qhop}--\eqref{eq:var_qhop} including \eqref{eq:rSEC2_qhop}, \eqref{eq:qhop_inflow}, \eqref{eq:layerConsecutive} and \eqref{eq:arcLastlayerLim}, \eqref{eq:nodeLastlayerLim} into our \scipjack framework. During the transformation into the layered graph, the following \textit{directed root cost test} decides if arcs are included in the layered graph (see, e.g., \textcite{gouveia2011}): If the arcs $(i,j), (r,j)\in A$ exist, $i \neq r$ and it holds that $c_{ij} \ge c_{rj}$, the arc $(i,j)$ is removed during the transformation, as node $j$ will always be connected to the root node. The shortest-path-based primal heuristic is extended to handle the layered graph, the hop-constraints, and the radial layout.

As before the bilinear products $X_{i i^{\prime}} X_{j j^{\prime}}$ in \eqref{eq:interf_tot_qhop} are linearized, but, instead of adding them all in the beginning as done before by \eqref{eq:quota_cons_trans_lin}~--~\eqref{eq:bilinear_var}, they are separated on-the-fly during the the B\&C algorithm. For the computational experiments in this section, we use \textit{pseudocost branching with strong branching initialization} presented in \textcite{Pedersen2025OR}. Until depth three of the B\&B tree the pseudcosts are initialized by performing strong branching on the ten fractional vertices which contribute the most interference in the current LP-solution. Afterwards, the pseudocosts are used. For more details, the reader is referred to \textcite{Pedersen2025OR, ACHTERBERG2005}. We call the resulting implementation \setting{QSTPI-HOP}.

According to \textcite{cazzaro2023}, each cable allows for a maximum of six 15\,MW wind turbines, resulting in a hop-limit $H=6$. In the flow-based formulation this cable capacity is simply introduced in Eq.~\eqref{eq:activeArc_flow}. We refer to this capacitated flow-based formulation as \setting{FLOW-CAPA}.

As a test set we consider the small- and medium-sized instances of Section~\ref{sec:data}. We use the same computational settings as before. We use the shifted geometric mean \parencite{achterberg2007Phd} with a shift of one second regarding time and a shift of 100 regarding the B\&B nodes. 

\begin{figure}
\centering
\includegraphics[width=.45\linewidth]{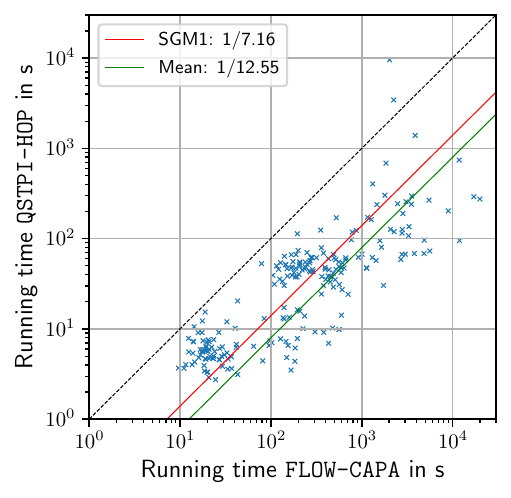}
\caption{Comparison of running times between \setting{QSTPI-HOP} and \setting{FLOW-CAPA} of all instances solved by both settings. The dotted line shows the break even line; the \textit{SGM1} line represents the shifted geometric mean of the speedup with a shift of one second; the \textit{Mean} line represents the arithmetic mean of the speedup.}
\label{fig:speedup-qstphop-grb}
\end{figure}

The results are summarized in Table~\ref{tab:hop_results}. The computational time for each individual instance needed by the two settings is shown in Fig.~\ref{fig:speedup-qstphop-grb}. Our proposed approach outperforms the flow-based formulation in all categories: It solves three more instances, it shows an average speedup of 7.16 in terms of time, and it reduces the number of B\&B nodes by a factor of 10.54. Furthermore, \setting{QSTPI-HOP} is on all but two instances faster than \setting{FLOW-CAPA}.

\begin{table}
    \centering
    \begin{tabular}{l|rrr}
    Setting & \# solved & time [s] & B\&B nodes\\
    \hline
    \texttt{QSTPI-HOP}& 223 &  27.20 &  85.11\\
    \texttt{FLOW-CAPA}& 220 & 194.77 & 897.29\\
    \texttt{QSTPI-Ilb}$^*$& 239 &   5.85 &  69.69\\
    \texttt{FLOW-Ilb}$^*$ & 228 & 384.08 & 2224.40
    \end{tabular}
    \caption{Computational results: The second column shows the number of instances solved to optimality; the third column shows the mean time over all instances solved by both settings for the capacitated approach using the shifted geometric mean with a shift of one second; Last column shows the number of average B\&B nodes over all instances solved by both settings using the shifted geometric mean with a shift of 100. *Row three and four summarize the results from the uncapacitated approach, see Section~\ref{sec:CompStudy}.}
    \label{tab:hop_results}
\end{table}

Comparing \setting{QSTPI-HOP} and \setting{FLOW-CAPA} to their uncapacitated versions, the flow-based approach benefits by introducing a capacity on the flow variables as the solution space is reduced. On the common solved instances, the computational time is reduced by around 50\%, although eight instances are solved less in total. The opposite can be observed for the \setting{QSTPI-HOP} as we solve an enlarged graph: The computing time is increased by a factor of around five compared to the \setting{QSTPI-Ilb}.

\section{Price of sequentiality}\label{sec:Res_Sequentiel_vs_combined}

The common approach of solving the wind farm layout and cable routing problem sequentially can lead to suboptimal solutions \parencite{fischetti2017Phd, cazzaro2022Phd, cazzaro2023}. As a final evaluation, we illustrate the price of that sequentiality. We compare:
\begin{itemize}
    \item the sequential approach of first finding a feasible positioning minimizing the turbine costs, followed by solving the capacitated routing problem (\setting{SEQ}),
    \item solving the capacitated routing problem on the positions found by the QSTPI (\setting{QSTPI-SEQ}), and
    \item solving the hop-constrained QSTPI (\setting{QSTPI-HOP}).
\end{itemize}
Concerning \setting{SEQ}, finding a feasible positioning $V^*$ minimizing the turbine costs is done by 
\begin{align}
\min\quad\sum_{i\in\potTerminals}\vertexCostsNumbered{i}\nodeVar_{i}\\
    \text{s.t.} \quad\eqref{eq:MinInterfQuota} - \eqref{eq:MinInterfEnd},
\end{align}
and is solved with \gurobiVersion{11}{01}. Having found $V^*$ by \setting{SEQ} or \setting{QSTPI-SEQ}, the capacitated routing problem is solved on the complete subgraph $G^*=(V^*, A^*, r)$ containing all chosen turbines $V^*$ and the substation $r$ and is formulated as
\begin{align}
    \min\quad &c^T x&\\
    \text{s.t.} \quad &\sum_{\arc\in\incomingArcs{\vertex}}\arcFlow_\arc - \sum_{\arc\in\outgoingArcs{\vertex}}\arcFlow_\arc = 1 & \forall\vertex\in V^*\setminus \{r\}\\
    &f_a\le Hx_a &\forall a\in A^*\\
    &x_a\in\lbrace0,1\rbrace, f_a\in\R_{\ge0}&\forall a\in A^*
\end{align}
and is solved with \gurobiVersion{11}{01}.

The test set consists of the 214 small- and medium-sized instances solved to optimality by both \setting{QSTPI-HOP} and \setting{QSTPI-Ilb}. The relative cost change is calculated by 
\begin{align*}
    c_{\mathrm{rel}}^{\mathrm{red}} = \left(1 - \frac{c_{\setting{SEQ}}}{c_{\setting{QSTPI}}}\right) \cdot100\%,
\end{align*}
where $c_\setting{SEQ}$ refers to the sequential approach of first minimizing turbine costs and the cable costs, and $c_\setting{QSTPI}$ to the respective \setting{QSTPI} approach.

\begin{figure}
    \centering
    \includegraphics[width=\linewidth]{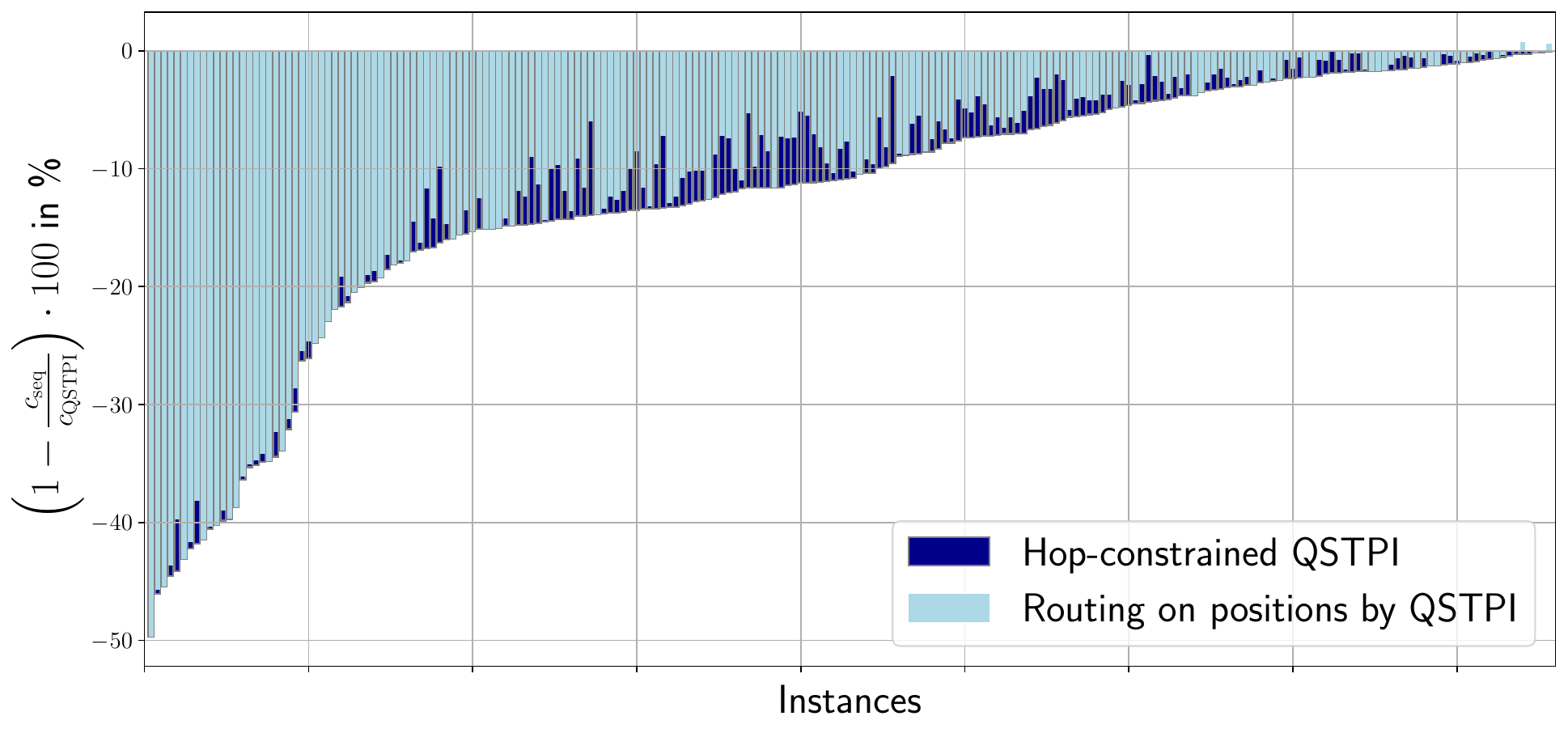}
    \caption{Cost reduction by using chosen turbines of \setting{QSTPI-SEQ} (light) and additionally by using \setting{QSTPI-HOP} (dark) compared to \setting{SEQ}}
    \label{fig:cost_red_seq_vs_seqQSTPI_vs_hopQSTPI}
\end{figure}

\begin{figure}
\centering
\begin{subfigure}[t]{.45\linewidth}
\includegraphics[width=.8\linewidth]{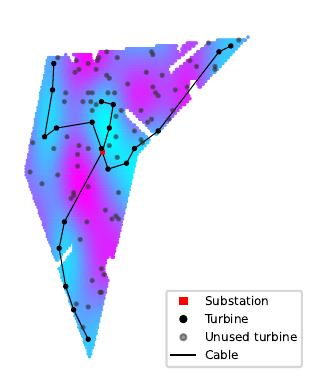}
\caption{Solution for \setting{SEQ}.}\label{fig:seq_example}
\end{subfigure}
\begin{subfigure}[t]{.45\linewidth}
\includegraphics[width=.8\linewidth]{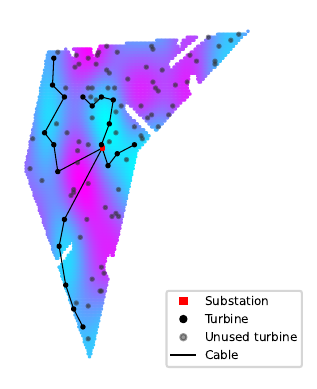}
\caption{Solution \setting{QSTPI-SEQ}; cost reduction 5.64\% compared to \setting{SEQ}.}\label{fig:seqQSTPI_example}
\end{subfigure}
\begin{subfigure}[b]{.45\linewidth}
\includegraphics[width=.8\linewidth]{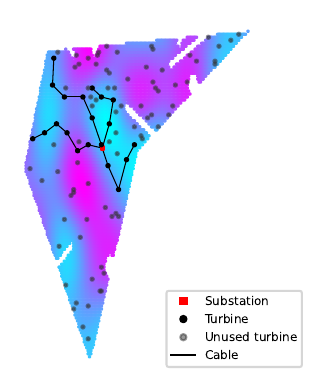}
\caption{Solution for \setting{QSTPI-HOP}; cost reduction 9.94\% compared to \setting{SEQ}.}\label{fig:hopQSTPI_example}
\end{subfigure}
\caption{Three solutions for the same hop-constrained QSTPI instance showing the price of sequentiality.}
\label{fig:seq-vs-integrated-example}
\end{figure}

Fig.~\ref{fig:cost_red_seq_vs_seqQSTPI_vs_hopQSTPI} presents the relative cost change of our novel approaches compared to the sequential approach for each instance. Already with \setting{QSTPI-SEQ}, the costs are reduced by 10.73\% on average, while two instances have a worse objective compared to the sequential approach, but with costs increasing by less than one percent. Using \setting{QSTPI-HOP} instead of the sequential approach the costs are reduced by 12.25\% on average, always improving the solution. In 43 instances the \setting{QSTPI-SEQ} gives the same solution as \setting{QSTPI-HOP}. Fig.~\ref{fig:seq-vs-integrated-example} shows the different solutions for an exemplary instance.

Already \textcite{Cazzaro.2022} showed that an integrated approach improves the solution quality of their heuristic. 
We conclude that the price of sequentiality is significant and justifies invoking an integrated approach. Therein, the QSTPI can provide both a lower bound on the optimal solution of the capacitated IWLCR problem as well as serve as a heuristic for the problem. 

\section{Conclusion and outlook}\label{sec:Discussion}
Current methods for solving the integrated wind farm layout and cable routing problem (IWFLCR) are either heuristic approaches or are limited to very small instances. This paper proposes a new approach for optimizing simultaneously the turbine locations and their connection to the grid's substation. The core is the QSTPI, an extension of the quota Steiner tree problem (QSTP). Besides the quota and connectivity constraints from the QSTP, the formulation includes the wake effect caused by a turbine reducing the profit of the surrounding turbines, and can also guarantee a minimum distance between chosen turbines. Both are vital technical constraints for the WFLO. We have upgraded the state-of-the-art STP-solver \scipjack so that it can handle the interference-constrained variant of the QSTP. Where the basic QSTP is of great advantage compared to standard flow-based formulation solvable by general MIP solvers, its interference-constrained variant bears no advantage at first. On the one hand this shows the immense difficulty of the quadratic constraints induced by the wake effect, and on the other hand the great improvements made by generic MIP solvers for these types of constraints. 

Inspired by known cutting planes for the interference constraint, we therefore propose a solution strategy of splitting the problem based on the total interference. Although two possibly challenging subproblems must be solved, we demonstrate on a large testset that a) the solution strategy speedups the solution process and makes most of the instances solvable, and b) the advanced version of \scipjack outperforms standard MIP solvers and is capable of finding optimal solutions for large-sized instances, where general state-of-the-art MIP solvers struggle. In particular, already using a small lower bound on the total interference results in a immense speed-up (see results for \setting{QSTPI-Ilb}). By splitting the problem using valid lower bounds calculated by an auxiliary MIP the number of solvable instances increases even further (see results for \setting{STPI-S-minI-1800}). Finally, we illustrate that to avoid excessive design costs for wind farm projects, an integrated view of the WFLO and WFCR is beneficial. This conclusion remains valid when cable capacity limits have to be respected.

Although our presented extensions to \scipjack outperforms general MIP solvers in terms of computational time as well as in solving large-sized instances, there are several instances of the considered dataset that have not been solved to optimality. As stated in Section~\ref{sec:implement}, the chosen primal heuristic often improves the primal solution only at the root node of the branch-and-bound tree, more efficient and suitable primal heuristics should be investigated. 

Furthermore, improving the dual bound of the presented formulation is vital. However, the cuts proposed by \textcite{fischetti2022} introduce very dense rows in the constraint matrix, making the LP-relaxation harder to solve, and preliminary experiments show that while improving the root node relaxation, the overall performance worsens. 

In this study, our novel approaches is validated on instances which are based on generic offshore wind farm areas which are characterized by a relatively uniform terrain -- although different foundation costs were considered --, wake expansion behavior, and using only a single turbine type. Nevertheless, we believe the proposed methodology is not strictly limited to the offshore case and could, for example, be applied for onshore wind farm planning in the future. By explicitly modeling heterogeneous terrains, varying turbine heights, complex wake calculations, and stricter layout constraints in the input data, the workflow of the QSTPI and its hop-constrained variant provides an interesting approach.

We believe that the hop-constrained QSTPI approach can be further sharpened to integrate further practically relevant constraints, e.g., substation capacity limits that limit the number of connected cable strings, or generalized to other types of network topologies.
Furthermore, it would be interesting to compare our optimization model with state-of-the-art simulation tools, e.g., FLORIS\footnote{\url{https://github.com/NREL/floris}, accessed on: 12/02/2025}.

{\singlespacing
\section*{Acknowledgements}
The work for this article has been conducted in the Research Campus MODAL funded by the German Federal Ministry of Education and Research (BMBF) (fund numbers 05M14ZAM, 05M20ZBM, 05M2025). We thank Davide Cazzaro for making his code available to compute the interference matrix for the used dataset.}
{\singlespacing 
\section*{CRediT statement}
Conceptualization: J.P., T.K.; data curation: J.P.; formal analysis: J.P., N.L.; investigation: J.P.; methodology: J.P., N.L., D.R.; software: J.P., D.R.; validation: J.P.; visualization: J.P.; writing – original draft: J.P.; writing – review and editing: J.P., N.L., D.R.}
{\singlespacing}

{\singlespacing
\printbibliography
}

\end{document}